\documentclass[11pt]{article}
\usepackage[utf8]{inputenc}

\usepackage[all, cmtip]{xy}
\usepackage{hyperref}
\usepackage{url}
\usepackage[utf8]{inputenc} % allow utf-8 input
\usepackage[T1]{fontenc}    % use 8bit T1 fonts
\usepackage{hyperref}       % hyperlinks
\usepackage{url}   
 \usepackage{parskip}
\usepackage{booktabs}       % professional-quality tables
\usepackage{amsfonts}       % blackboard math symbols
\usepackage{nicefrac}       % compact symbols for 1/2, etc.
\usepackage{microtype}      % microtypography
\usepackage[pdftex]{graphicx}
\usepackage{todonotes}
\usepackage{amsmath}

\usepackage{tikz-cd}
\usepackage{amsfonts}
\usepackage{amssymb}
\usepackage{amsthm}
\usepackage{hyperref}
\usepackage[noadjust]{cite}
\usepackage{caption}
\usepackage{bbm}
\usepackage[scr=boondoxo]{mathalpha}

\usepackage{tikz}
\usetikzlibrary{decorations.pathreplacing}
\usepackage[T1]{fontenc}
\usepackage[utf8]{inputenc}
\usepackage{mathtools}
\usepackage{verbatim}

\usepackage[top=2.9cm, bottom=2.9cm, left=2.5cm, right=2.5cm, headsep=0.2in]{geometry}

\usepackage{fancyhdr}
\usepackage{cleveref}
\title{Equivariant Intrinsic  formality}
\author{ST and RS }
\theoremstyle{plain}
\newtheorem{theorem}{Theorem}[section]

\newtheorem{corollary}[theorem]{Corollary}

\newtheorem{prop}[theorem]{Proposition}

\newtheorem{lemma}[theorem]{Lemma}

\theoremstyle{definition}
\newtheorem{definition}[theorem]{Definition}

\newtheorem{exmp}[theorem]{Example}
\theoremstyle{remark}
\newtheorem{cons}[theorem]{Construction}
\newtheorem{remark}[theorem]{Remark}

\newcommand{\QQ} {\mathbb Q}

\def\sO{\mathscr{O}}
\def\sU{\mathscr{U}}
\def\sV{\mathscr{V}}
\def\sB{\mathscr{B}}
\def\sN{\mathscr{N}}
\def\sM{\mathscr{M}}
\def\sA{\mathscr{A}}
\def\sE{\mathscr{E}}
\def\sI{\mathscr{I}}
\def\cA{\mathcal{A}}
\def\cM{\mathcal{M}}
\def\cB{\mathcal{B}}
\def\cV{\mathcal{V}}
\def\cI{\mathcal{I}}

\def\cU{\mathcal{U}}

\def\C{\text{C}}
\def\PH{\text{PH}}
\def\sP{\mathscr{P}}
\def\sK{\mathscr{K}}
\def\coker{\text{coker}}
\def\S{\text{S}}
\def\ker{\text{ker}}
\def\DGA{\text{DGA}}

\newcommand{\beq} {\begin{equation}}
\newcommand{\eeq} {\end{equation}}

\renewcommand{\Im} {\operatorname{Im}}

\begin{document}

	\title{Equivariant Intrinsic Formality}

	\author{ Rekha Santhanam\footnote{Department of Mathematics, Indian Institute of Technology Bombay, India. reksan@iitb.ac.in}, Soumyadip Thandar\footnote{Department of Mathematics, Indian Institute of Technology Bombay, India.
 stsoumyadip@gmail.com} }
	\maketitle
	%\maketitle 

    \begin{abstract}

Algebraic models for equivariant rational homotopy theory were developed by Triantafillou
and Scull
for finite group actions and $S^
1$
action, respectively. They showed that given a diagram of rational
cohomology algebras from the orbit category of a group $G$, there is a unique minimal system of $\DGA$s
representing a unique $G$-rational homotopy type that is weakly equivalent to it. However, there can be
several equivariant rational homotopy types with the same diagram of cohomology algebras. Halperin,
Stasheff, and others studied the problem of classifying rational homotopy types up to cohomology in the
nonequivariant case. In this article, we consider this question in the equivariant case. For the case $G=C_p$,
for prime $p$, under suitable conditions, we are able to determine the equivariant rational homotopy types
with isomorphic diagram of cohomology algebras in terms of non-equivariant data. We give explicit
examples to demonstrate how these theorems can be applied to classify equivariant rational homotopy
types with isomorphic cohomology.

%Algebraic models for equivariant rational homotopy theory were developed by Triantafillou \cite{GT82} and Scull \cite{LS02}, \cite{LS08} for finite group actions and $S^1$ action, respectively.  They showed that given a diagram of rational cohomology algebras from the orbit category of a group $G$, there is a unique minimal system of $\DGA$s representing a unique $G$-rational homotopy type that is weakly equivalent to it. However, there can be several equivariant rational homotopy types with the same diagram of cohomology algebras. Halperin, Stasheff, and others (\cite{HARST}, \cite{SY82}, \cite{LUPT}, \cite{HIRO}, \cite{SS12}) studied the problem of classifying rational homotopy types up to cohomology in the non-equivariant case. In this article, we consider this question in the equivariant case. For the case $G=C_p$, for prime $p$, under suitable conditions, we are able to determine the equivariant rational homotopy types with isomorphic diagram of cohomology algebras in terms of non-equivariant data. We give explicit examples to demonstrate how these theorems can be applied to classify equivariant rational homotopy types with isomorphic cohomology.

	\end{abstract}

\maketitle

\section{Introduction}

 Any two simply connected spaces are said to have the same \emph{rational homotopy type} if there is a zigzag of morphisms between them, each inducing isomorphism on their rational cohomology. Quillen and Sullivan give algebraic models, namely,  differential graded Lie algebras \cite{QR} and differential graded commutative algebras, written as $\DGA$s in short (minimal algebras, \cite[Section 2]{Sul77}) respectively, describing simply connected spaces up to their rational homotopy type.

We say two $G$-simply connected spaces (i.e., $G$-space $X$ whose fixed point spaces, $X^H$ are simply connected for all subgroups $H$ of $G$) have the same \emph{ $G$-rational homotopy type} if there is a zigzag of $G$-maps, each inducing an isomorphism on the rational cohomology of the fixed point  of the spaces under every subgroup of $G$.
 Triantafillou \cite{GT82} (for finite $G$) and Scull \cite{LS02} (for  $G=\S^1$) define algebraic models describing $G$-simply connected spaces up to the same  $G$-rational homotopy type. These models lie in the subcategory of injective objects  (\Cref{defn:inj object}) of the category of functors from the orbit category of $G$, $\sO_G$  (\Cref{defn:orbit}), to the category of cohomologically  $1$-connected $\DGA$s. We refer to functors  from $\sO_G$, to the category of cohomologically  $1$-connected $\DGA$s (graded algebras/ vector spaces) as a \emph{diagram of  $\DGA$s} (graded algebras/ vector spaces) over $\sO_G$  and as a \emph{system of $\DGA$s}  over $\sO_G$ when it is injective. In \cite[Theorem 6.2]{GT82}, \cite[Theorem 6.28]{LS02}, the authors show that there is a one-to-one correspondence between  $G$-simply connected spaces up to the $G$-rational homotopy type and isomorphism classes of minimal system of $\DGA$s (\Cref{def:minimal}) over $\sO_G$.  

A natural question is to ask for  all minimal algebras and, therefore, rational homotopy types with isomorphic (rational) cohomology algebras. Given a graded algebra  $A^{\ast}$, define the moduli set of all minimal algebras (up to isomorphism) with cohomology $A^{\ast}$, $$\cM_{A^\ast}:=\{\text{rational homotopy type of X}~|~H^\ast(X;\QQ)\cong A^\ast\}.$$
In the non-equivariant context, the  set $\mathcal{M}_{A^\ast}$ has been studied by several authors including \cite{HARST}, \cite{SY82}, \cite{LUPT}, \cite{HIRO}, \cite{SS12} with different view points. Lupton \cite{LUPT} shows that for any positive integer $n$ there is a graded algebra $A^\ast$ such that the cardinality of $\mathcal{M}_{A^\ast}$ is $n$. 

Lemaire and Sigrist \cite{LS78} produce an infinite family of distinct rational homotopy types with the same cohomology algebra and rational homotopy Lie algebra.
Halperin and Stasheff \cite{HARST}, study $\mathcal{M}_{A^\ast}$ by considering the set of perturbations of a bigraded model constructed from $A^{\ast}$. In particular, they show  that for $A^\ast=H^\ast((S^2\vee S^2)\times S^3;\QQ)$, the set $\mathcal{M}_{A^\ast}$ consists of two points. Shiga and Yamaguchi \cite{HIRO} study the set $\mathcal{M}_{A^{\ast}}$ by constructing a correspondence between $\cM_{A^{\ast}}$ and rational points of Grassmann manifolds modulo an equivalence relation generated by the group of automorphisms of $\DGA$s (\cite[Theorem 2.1, Corollary 2.3]{HIRO}).

In this article, we study the equivariant analogue of $\cM_{A^\ast}$, defined as the moduli set of all minimal systems of $\DGA$s (up to isomorphism) with cohomology diagram of graded algebras $\cA^{\ast}$ over $\sO_G$,
$$\cM^{G}_{\cA^\ast}:=\{\sM \ | \ \sM \text{ is a minimal system of DGAs over}~\sO_G~ \text{and} ~H^\ast(\sM;\underline{\QQ})\cong \cA^\ast\},$$ 
where $\underline{\QQ}$ is the constant coefficient system defined by $\underline{\QQ}(G/H):=\QQ$, for every subgroup $H$ of $G$.
A minimal system of $\DGA$s over $\sO_G$ which determines the $G$-rational homotopy type of a $G$-simply connected space is obtained by taking elementary extensions  (\Cref{def:ele}) inductively. An elementary extension is the  equivariant analogue of the Hirsch extension used to construct minimal algebras (\cite[Chapter 16.2]{PGJM}).  
Unlike the non-equivariant case, the generators added at $n$-th stage extension in the construction of $\sM$, which we denote by $\sM_n$, can have degree greater than $n$.

 The construction of an elementary extension of a system of $\DGA$s $\cU$ over $\sO_G$, depends on the following data; a diagram of vector spaces $\underbar{V}$ over $\sO_G$ of degree $n$ and an element $[\alpha]\in H^n(\cU;\underbar{V})$, and the extension is denoted by $\cU^{\alpha}(\underbar{V})$. 
  Any two non-isomorphic minimal systems of $\DGA$s, with an isomorphic diagram of cohomology algebras, differ at some $n$-th stage. 

A necessary condition for isomorphic elementary extensions over the same system of $\DGA$s over $\sO_G$, is given by Scull  \cite[Proposition 11.52]{LS02}.

In \Cref{section:intrinc}, we define  \emph{Condition $\C_n$} on a system of $\DGA$s  which ensures that two elementary extensions over the same system of $\DGA$s are non-isomorphic, as proved in \Cref{prop:plural}. This improves our  understanding of $\cM^G_{\cA^{\ast}}$ and more specifically gives a method to construct  minimal systems of $\DGA$s that are not quasi-isomorphic but have the isomorphic cohomology diagrams. For instance, we  show that  $\cM^G_{\cA^{\ast}}$ can have more than  one point in  \Cref{example:plural}.

In order to study $\cM^G_{\cA^\ast}$, for a cohomology diagram of graded algebras $\cA^\ast$ over $\sO_G$ which can be considered as a  diagram of $\DGA$s with zero differential, it is imperative to understand its  minimal model, that is, a  minimal system of $\DGA$s $\sM$ over $\sO_{G}$ with a morphism $\rho:\sM\to \cA^{\ast}$ inducing isomorphism in cohomology. 
However, the cohomology diagram of a given $G$-space is not always injective, as can be seen from examples (\Cref{exmp:noninjective}). In such cases, we need to consider, the injective envelope  of the diagram of $\DGA$s (\Cref{thm:fineee}), whose differential need not be zero, making the ensuing computations more complex.  In \Cref{prop:injcohom}, we give a simple condition to verify when a diagram of $\DGA$s over an orbit category of $C_p$, for prime $p$, is injective.

Moreover, we observe that at each stage of the  construction of a minimal system if the \emph{associated} diagram of vector spaces (\Cref{def:associ:vect}) involved in the elementary extension is injective, then 
the computations simplify. We therefore, consider the following question:

\emph{What conditions on an injective cohomology diagram of graded algebras ensure that the associated diagram of vector spaces added for each elementary extension are injective?}
    
We \emph{answer} this question when $G=C_p$, where $p$ is prime. In \Cref{Z_P injimpliesminimal}, we show that if the cohomology diagram $\cA^\ast$ over $\sO_{C_p}$ has the property that the structure map $\cA(G/e)\to \cA(G/G)$ is a retract (\Cref{defn: retract}) and the minimal model is $\sM$, then the associated diagram of vector spaces at every stage of the construction of  $\sM$ is injective. In this case, the equivariant minimal model $\sM$ of $\cA^\ast$ is level-wise minimal, i.e., $\sM(G/H)$ is a minimal model for $\cA^{\ast}(G/H)$, reducing the problem to the non-equivariant case.

The diagram of graded algebras, $\cA^{\ast}$ over $\sO_G$, is said to be \emph{equivariantly $k$-intrinsically formal} if  there exists a $k$-isomorphism (\Cref{defn:n-iso}) between any two minimal systems in the equivariant moduli set $\cM^{G}_{\cA^\ast}$. Further, if  this set $\cM^{G}_{\cA^\ast}$ is a singleton, then we say $\cA^{\ast}$ is \emph{equivariantly intrinsically formal}.

Define a system of $\DGA$s to be \emph{equivariantly formal} if its weak equivalence class (\Cref{def:weakyeqivalent}) can be completely determined by its cohomology diagram. A $G$-space is equivariantly formal if the minimal system of $\DGA$s corresponding to it is equivariantly formal. With this set up, we consider the following question:

\emph{Can we compute the cardinality of $G$-rational homotopy types with isomorphic cohomology diagrams over $\sO_G$ or say when the cohomology diagram is equivariantly intrinsically formal?  }

We \emph{address} this  question for $G=C_p$, for $p$ prime, in \Cref{section:minimalmodels}. We extend the results of \cite{HIRO}, for systems of $\DGA$s over $\sO_{C_p}$ in Theorem \ref{theo:MAIN}.
 In \Cref{thm:cardofM}, under suitable conditions, we determine the cardinality of a subclass of $\mathcal{M}_{\cA^{\ast}}^{C_p}$  in terms of the non-equivariant set $\mathcal{M}_{\cA^{\ast}(C_p/e)}$. In \Cref{cor:main}, we give a sufficient condition for a diagram of graded algebras to be equivariantly intrinsically formal, that is, the equivariant moduli set corresponding to $\cA^{\ast}$ is a point. More precisely, we prove the following.

\textbf{Theorem} \ref{thm:cardofM} and \textbf{Corollary} \ref{cor:main}.
Let $\cA^{\ast}$ be a diagram of graded algebras over $\sO_{C_p}$ such that  its structure map is a retract  with minimal model $\sM$. Assume that $\cA^\ast$ is  equivariantly $(n-1)$-intrinsically formal and $\sM(C_p/C_p)$ does not have elements of degree $\geq n$. Then the following statements are true. 
  \begin{enumerate}
\item The set of isomorphism classes of minimal systems containing $\sM_{n-1}$ is determined by the moduli set corresponding to $\cA^{\ast}(C_p/e)$. 
  \item If $\cA^{\ast}(C_p/e)$ does not have elements of degree $>n+1$, then the cardinality of the equivariant moduli set corresponding to $\cA^{\ast}$ coincides with that of the moduli set corresponding to $\cA^{\ast}(C_p/e)$.

  \item If the minimal models of $\sM(C_p/e)$ and $H^{\ast}(\sM(C_p/e))$ are isomorphic, then $\cA^{\ast}$ is equivariantly intrinsically formal. In particular, $\sM$ is equivariantly formal.
 \end{enumerate}

This allows us to produce examples of equivariantly formal $C_p$-spaces, 
for instance, see \Cref{example sec5}. As a further application of \Cref{thm:cardofM} and \Cref{cor:main},  in \Cref{oddwedgecross},  we demonstrate how to prove that a given diagram of $\DGA$s is equivariantly $n$-intrinsically formal for a particular $n$.  Further,  using  the work of Shiga and Yamaguchi \cite{HARST} for the non-equivariant case, we are able to compute the cardinality of $\cM^{C_p}_{\cA^\ast}$ for a given diagram of graded algebras $\cA^{\ast}$ over $\sO_{C_p}$ in  \Cref{exmp: 3} and \Cref{exmp:last}.

In the forthcoming works \cite{SST},\cite{ST25}, we extend these results for $C_{p^n}$, $C_{pq}$ and $C_p\oplus C_p$-diagrams of graded algebras and compute other classes of equivariantly formal spaces.

\section{Background}\label{section:background}

In this article, we work with $\DGA$s over $\QQ$ and assume   $G$ to be a finite group. 
\begin{definition}\label{defn:orbit}
Given a group $G$, the category of canonical orbits is the category whose objects are $G$-sets $G/H$ and morphisms are $G$ maps between them. We denote this category by $\sO_G$. 
\end{definition}

\begin{definition}\label{defn:inj object}
    
An object $I$ in a category $\mathcal{C}$ is said to be \textit{injective} if for every injective morphism $f:X\to Y$ and every morphism $g:X\to I$

 \[
\xymatrix{
X \ar[d]_g\ar[r]^f& Y\ar@{-->}[dl]^{h} \\
I } 
\]

there exists a morphism $h:Y\to I$ such that $h\circ f=g$. 
\end{definition}

A {\it  diagram of $\DGA$s} is a  covariant functor from the orbit category $\sO_G$ to the category of cohomologically $1$-connected $\DGA$s. If this functor is injective, then we refer to it as a {\it system of $\DGA$s} in line with  \cite{GT82}. We will denote the category of  systems of $\DGA$s by  $\DGA^{\sO_G}$.

By forgetting the differential in a diagram of $\DGA$s we get a diagram of rational vector spaces, also known as a  {\it dual rational coefficient system}. The category of dual rational coefficient systems will be denoted by $Vec^\ast_G$.  
A {\it rational coefficient system} is a  contravariant functor from $\sO_G$ to the category of  rational vector spaces. We denote  the category of rational coefficient systems by $Vec_G$.\\

Let $X$ be a $G$-space such that $X^H$ is non-empty and simply connected for all $H\leq G$. Then the corresponding  diagram of cohomology algebra of $X$  with differential $0$ is $1$-connected (i.e., $H^1(X^H;\QQ)=0$ for every subgroup $H$ of $G$). This need not be an injective dual coefficient system. However, every dual coefficient system has an injective envelope.

We now describe \cite[Prop. 7.34]{LS02}, the  embedding of a given coefficient system  $\cA$ into its injective envelope $\cI$.

\begin{definition}\label{equation:24}

We define 
\beq 
V_H:=\cap_{H\subset K}\ker \cA(\hat{e_{H,K}}),
\eeq

 where $\hat{e_{H,K}}:G/H\to G/K$ is the projection and $\cA(\hat{e_{H,K}})$ is the induced structure map on the functor $\cA$.  Note that  $V_G$ is defined to be $\cA(G/G)$. Let $\cI=\oplus_H \underbar{V}^{\ast}_H$, where 
 \beq
 \underline{V}_H^{\ast}(G/K):=Hom_{\QQ(N(H)/H)}((G/H)^K,V_H).
 \eeq

There is an injective morphism $\cA\to \cI$ extending the natural inclusions of $\cap_ {H\subset K}\ker \cA(\hat{e_{H,K}})$.\\

\end{definition}
\begin{prop}\cite[Section 4]{GT82}
    A dual coefficient system $\cA$ is injective if and only if it is of the form $\cA = \bigoplus_H \underline{V}^{\ast}_H$ for some collection of $\QQ(N(H)/H)$-modules
$V_H$ and $$\underline{V}_H^{\ast}(G/K)=Hom_{\QQ(N(H)/H)}((G/H)^K,V_H).$$
\end{prop}

Given a diagram of $\DGA$s, forgetting the differential will give a dual rational coefficient system whose injective envelope is a diagram of $\DGA$s, with $0$ differential. However, the map into the injective envelope of dual rational coefficient system will not be a quasi-isomorphism in general. % 

Fine and Triantafillou \cite{FINE}, prove the existence of \emph{injective envelope} for a diagram of $\DGA$s.

\begin{theorem}\cite[Theorem 1]{FINE}\label{thm:fineee}
For a diagram of $\DGA$s $\sA$ over $\sO_G$, where $G$ is finite group, there is an injective system  of $\DGA$s  $\sI(\sA)$, called the injective envelope of $\sA$ along with an inclusion $i:\sA\to \sI(\sA)$ which is a quasi-isomorphism. 
\end{theorem} 

We now describe their construction.

  \begin{definition}\label{def:associ}
Let $G$ be a group and $H\leq G$. Let $A_H$ be a $\DGA$ over $\QQ$ such that $N(H)/H$ acts on it by $\DGA$ automorphisms. The \emph{associated system of $\DGA$s} $\sA_H$, of the $\DGA$ $A_H$, is a system of $\DGA$s defined as follows:
Let $V_H$ be a copy of $A_H$ considered as a graded $\QQ(N(H)/H)$-module by forgetting the differential and let $\underbar{V}_H^{\ast}$ be the induced injective diagram of vector spaces \Cref{equation:24}. Let $s\underbar{V}_H^{\ast}$ be copy of $\underbar{V}_H^{\ast}$ with a shift of degree by $+1$. We denote by $\bigwedge_H$ the system of acyclic $\DGA$s generated by $\underbar{V}_H^{\ast}\oplus s\underbar{V}^{\ast}_H$, where $d(\underbar{V}^{\ast}_H)=s\underbar{V}^{\ast}_H$. Now we define the associated system $\sA_H$ by

\[
\sA_H(G/K)=
\begin{cases}
     \bigwedge_H(G/K), & \text{for } (K)<(H)\\
    Hom_{\QQ(N(H)/H)}(\QQ(G/H)^K, A_H),  & \text{ for } (K) \neq (H)
\end{cases}
\]

where $(H)$ is the conjugacy class of $H$ in $G$. The value of this functor on morphism is obvious.

\end{definition}

\begin{definition}
Let $\sA$ be a system of $\DGA$s and let $A_H$ be the subalgebra of $\sA(G/H)$ which is equal to $\cap _{H\subset H'}\ker \sA_{H,H'}$, where $\sA_{H,H'}$ is the morphism induced by the projection $G/H\to G/H'$. Let $\sA_H$ be the associated system to $A_H$. The enlargement of $\sA$ at $H$ is the system of $\DGA$s $\sI_H(\sA)$ defined by 

\[
\sI_H(\sA)=
\begin{cases}
     \sA(G/K)\otimes \sA_H(G/K), & \text{for } K<H\\
    \sA(G/K),  & \text{otherwise } 
\end{cases}
\]

where $K<H$ means that $K$ is a proper subgroup of a conjugate of $H$. The value of the functor $\sI_H(\sA)$ on morphisms is the obvious one, namely, they are equal to the old morphisms when restricted to the subsystems $\sA$ and $\sA_H$ respectively.

The injective envelope of a system of $\DGA$s $\sA$ is constructed by the following steps. First, we consider enlargement at $G$, namely, $\sI_G(\sA)$ of the given system. Next we  consider a maximal subgroup $H$ of $G$, and take the enlargement of $\sI_G(\sA)$ at $H'$, (that is, we construct $\sI_{H'}(\sI_G(\sA))$,) where $H'$ is some conjugate of $H$. We repeat this process untill we reach the trivial subgroup. For details, see \cite{FINE}. 

\end{definition}

\begin{remark}
 Note in the construction of the injective envelope, we add new elements and kill their cohomology class by adding their suspension. So if we start with an injective diagram of graded algebras, thought of as a diagram of $\DGA$s by considering the differential zero, then it is injective as a diagram of $\DGA$s.

\end{remark}

We now want to  define weak equivalences on the category  $\sO_G[\DGA]$. Before that, we define  the notion of homotopy on systems of $\DGA$s.

\begin{definition}

Given a system  of $\DGA$s $\sU$, define $\sU(t,dt)$ as the $\DGA$ diagram $$\sU(t,dt)(G/H)=\sU(G/H) \otimes_\QQ \QQ(t,dt), $$ where $\QQ(t,dt)$ are free with $t$ in degree $0$ and $dt$ in degree $1$. Two morphisms of $\DGA$s $f,g: \sU_1 \to  \sU_2$ are  said to be homotopic if there exists a $\DGA$ morphism 
 $H:\sU_1\to \sU_2\otimes \mathbb{Q}(t,dt)$ such that $p_0H=f$ and $p_1H=g$ where $p_i:\sU_2\otimes \mathbb{Q}(t,dt)\to \sU_2$ are defined as $p_i(t)=i$ for $i=0,1$ and $p_i(dt)=0$.

\end{definition}
This does not define an equivalence relation on $\DGA^{\sO_G}$. However, it does give an equivalence relation on {\it minimal} (Definition \ref{def:minimal}) systems of $\DGA$s, which we discuss later in this section. We define a coarser relation on $\DGA^{\sO_G}$. Given two diagrams of $\DGA$s $\sU, \sB$,  if there is a morphism  $f$ in $\sU \to \sB$ or $\sB \to \sU$ inducing a cohomology isomorphism at each level (at $G/H$ for all $H \leq G$) then $f$ is said to be a quasi-isomorphism. The equivalence relation generated by quasi-isomorphisms is defined as a weak equivalence of systems of $\DGA$s. Using the notion of injective envelopes, we can define weak equivalence on the category of diagrams of $\DGA$s

\begin{definition}\label{def:weakyeqivalent}
Let $\cU$ and $\cV$ be two $\DGA$ diagrams over $\sO_G$. We say $\cU$ and $\cV$ are  \emph{weakly equivalent} if there is a weak equivalence between their injective envelopes.
\end{definition}

Recall that associated with any $G$-space $X$, there is the system of $\DGA$s given by the de Rham-Alexander-Spanier algebra $\sE(X)(G/H):=\cA(X^H)$ for every $H\leq G$. Triantafillou \cite[Theorem 1.5]{GT82} proves that there is a bijective correspondence between the $G$-space $X$ (with every fixed point set simply connected) and the minimal system of $\DGA$s $\sM_X$ of $\sE(X).$

Scull generalizes these ideas to spaces with an $S^1$ action. In \cite[Section 21]{LS02}, Scull shows that, unlike the non-equivariant case, the notion of minimality in the equivariant case arising from filtration via minimal extensions of systems of $\DGA$s does not satisfy the decomposability condition. 

Note that homotopy defines an equivalence relation on morphisms from $\sM \to \sB$ for any system of $\DGA$s $\sB$, whenever $\sM$ is a minimal system \cite[Prop. 3.5]{LS02}. Further, given a quasi isomorphism $\rho: \sU \to \sB$ of a system of $\DGA$s and a morphism $f:\sM\to \sB$ from a minimal system $\sM$, there is a lift $g:\sM\to \sU$ such that $\rho g\simeq f$, \cite[Prop. 3.6]{LS02}.

The results (\cite[Theorem 1.5]{GT82}, \cite[Theorem 4.13]{LS02}) show that the category of $G$-spaces (whose fixed points sets are simply connected)  up to rational homotopy equivalences is equivalent to the category of minimal systems of $1$-connected $\DGA$s modulo homotopy equivalences.

In order to give the construction of a minimal model of a system of $\DGA$s we first define elementary extensions.  

\begin{definition}\label{def:ele}
Given a system of $\DGA$s $\mathcal{U}$, a diagram of vector spaces $\underline{V}$ assigned to be of degree $n$, and a map $\alpha: \underline{V}\to \underline{Z}^{n+1}(\mathcal{U})$ (here $\underline{
Z}(\cU)$ denotes the kernel of $\cU$), the {\it elementary extension } of $\mathcal{U}$ with respect to $\alpha$ and $\underbar{V}$, denoted by $\mathcal{U}^{\alpha}(\underbar{V})$, is constructed as follows.

Let $\underbar{V}\to \underbar{V}_0 \xrightarrow{w_0} \underbar{V}_1\xrightarrow{w_1} \underbar{V}_2\cdots$
be a minimal injective resolution of $\underbar{V}$ constructed by taking $\underbar{V}_i$ to be the injective embedding of $\coker w_{i-1}$, which is of finite length.  

Construct a commutative diagram 

\[
\xymatrix{
\underline{V}\ar[d]_{\alpha} \ar[r] & \underline{V}_0\ar[d]^{ \alpha_0}\ar[r]^{w_0}& \underline{V}_1\ar[d]^{\alpha_1}\ar[r]^{w_1} & \underline{V}_2\ar[r]\ar[d]^{\alpha_2} & \cdots  \\
\underline{Z}^{n+1}(\mathcal{U})\ar[r]         & \mathcal{U}^{n+1}\ar[r]_{d}&\mathcal{U}^{n+2}\ar[r]_d & \mathcal{U}^{n+3}_d \ar[r] &\cdots
}
\] 

The  maps $\alpha_i$ are  constructed inductively by  first noting that $d\alpha_iw_{i-1}=dd\alpha_{i-1}=0$, so $d\alpha_i|_{\Im w_{i-1}}=0$ and then by the injectivity of $\cU$ we get a commutative diagram:

\[
 \xymatrix{
\underbar{V}_i/\Im w_{i-1} \ar[d]_{d\alpha_i}\ar@{^{(}->}[r]^(0.6){ \rho^{\ast}} &\underbar{V}_{i+1}\ar[dl]^{\alpha_{i+1}}\\
\cU^{n+i+1} 
}
\]

Define $\cU^{\alpha}(\underbar{V}):=\cU\otimes(\otimes_i\QQ(\underbar{V}_i))$, where $\QQ(\underbar{V}_i)$ is the free graded commutative algebra generated at $G/H$ by the vector space $\underbar{V}_i(G/H)$ in degree $n+i$; the differential is defined on $\cU$ by the original differential on $\cU$, and on the generators of $\underbar{V}_i$ by $d=(-1)^i\alpha_i+w_i$. Since  $\underbar{V}_i$ is injective for all $i$ by construction, as a vector space the system is the tensor product of injectives and hence injective. Thus, $\cU^{\alpha}(\underbar{V})$ is a new system of $\DGA$s.

\end{definition}

\begin{remark}
We use $\cU^{\alpha}(\underbar{V})$ to denote the elementary extension of $\cU$ by a diagram of vector spaces $\underbar{V}$ with respect to the map $\alpha$. We denote the elementary extension of $\cU$ by the diagram of vector spaces $\underbar{V}$ by $\cU(\underbar{V})$ if we do not want to focus on the map $\alpha$.
\end{remark}

The following result shows when two elementary extensions are isomorphic. 

 \begin{prop}\cite [Lemma 11.53]{LS02}\label{prop:11.53}
  Suppose $f:\mathcal{U}^{\alpha}(\underline{V})\to \mathcal{U}^{\alpha'}(\underline{V'})$ is a map between two degree $n$ elementary extensions of $\mathcal{U}$ with the following properties: 
 \begin{enumerate}
  \item $f$ restricts to an isomorphism of $\mathcal{U}$.
      \item On $\underline{V}$, $f(x)=g(x)+a(x)$, where $g:\underline{V}\to \underline{V'}$ is an isomorphism and $a(x)\in \mathcal{U}$. 
       \end{enumerate}
 Then $f$ is an isomorphism.  
  \end{prop}

  A minimal system of $\DGA$s is  defined as follows. 
\begin{definition}\label{def:minimal}
A system of $\DGA$s $\sM$ is minimal if $\sM=\cup_n\sM_{n}$, where $\sM_{0}=\sM_{1}=\underline{\mathbb{Q}}$ and $\sM_{n}=\sM_{n-1}(\underline{V})$ is the elementary extension for some diagram of vector spaces $\underline{V}$ of degree $\geq n$. 
\end{definition}

\begin{theorem}\label{thm:quasiisiso}[\cite{LS02} Theorem 3.8]
If $f:\sM\to \sN$ be a quasi-isomorphism between minimal systems of $\DGA$s, then $f\simeq g$, when $g$ is an isomorphism. 

\end{theorem}

Thus, if we have two minimal systems $\sM$, $\sN$ and quasi-isomorphisms $\rho_1:\sM\to \cU$ and $\rho_2:\sN\to \cU$ by the lifting property of maps from minimal systems to systems of $\DGA$s we get a map $f:\sM\to \sN$ which is a quasi-isomorphism. By Theorem \ref{thm:quasiisiso}, we get $f\simeq g$ where $g$ is an isomorphism. Now we define the following.

\begin{definition}\label{def:minimalmdel}
If $\sM$ is a minimal system and $\rho: \sM\to \sU$ is a quasi-isomorphism, we say that $\sM$ is a minimal model of $\sU$. 
\end{definition}

Maps between two minimal systems of $\DGA$s are much {\it nicer}, in the sense that they are always homotopic to a level-wise map of extensions. We will make use of this fact later. 
\begin{lemma}\cite [Lemma 13.57]{LS02}\label{lemm:13.57}
Any morphism $f:\sM\to \sN$ between minimal systems of $\DGA$ is homotopic to a morphism $g$ which maps $\sM_{n}$ to $\sN_n$ for all $n$.

\end{lemma}
\begin{remark}\label{rmk:isolevelwise}
Given a morphism $g$ as in \Cref{lemm:13.57}, by \Cref{thm:quasiisiso}, we get that $g$ is an isomorphism. Then \Cref{prop:11.53} implies that if $\sM_{n}=\sM_{n-1}(\underline{V})$ and $\sN_n=\sN_{n-1}(\underline{V'})$ then $\underline{V}\cong \underline{V'}$ and for any $x \in \underline{V}$, $g(x)=x +b$ for some $b \in \sN_{n-1}$.    
 \end{remark}

Observe that any minimal system is cohomologically $1$-connected, that is, it satisfies $\underline{H}^0(\sM)=\mathbb{Q}$ and $\underline{H}^1(\sM)=0$. It can be shown that being cohomologically $1$-connected is sufficient for a diagram of $\DGA$s to have a {\it minimal model.}

\begin{theorem}\cite [Theorem 3.11]{LS02}
If $\sU$ is a system of $\DGA$s which is cohomologically $1-$connected, then there exists a minimal model of $\sU$, i.e., a minimal system $\sM$ and a quasi-isomorphism $\rho:\sM \to \sU$.
\end{theorem}

Note that this construction ensures that $\rho$ restricted to $\sM_n \to \sU$ is an \emph{$n$-isomorphism}. 

\begin{definition}\label{defn:n-iso}
    We say a morphism $f:\cU \to \cB$ between two systems of $\DGA$s is an $n$-isomorphism if $f^{\ast}:H^{\ast}(\cU)\to H^{\ast}(\cB)$ is isomorphism up to degree $n$ and monomorphism at degree $(n+1)$.  
\end{definition}

In the non-equivariant setup, we define the following.
\begin{definition}[Retract]\label{defn: retract}
   Given a $\DGA$  $A$ and a sub-$\DGA$ $B$ of $A$, we say $B$ is a retract of $A$ if there is a $\DGA$-morphism $r: A\to B$ such that $r\circ i=id_{B}$. Here $i:B\to A$ is the inclusion morphism and the morphism $r$ is called the retraction. 
   %Equivalently, we define a map $r:A \to B$ is said to be a retraction if there is a $\DGA$ morphism $i:B\to A$ such that $r\circ i=id_{B}$. 
\end{definition}

\section{ Equivariantly intrinsically formal graded algebras} \label{section:intrinc}
A cohomologically $1$-connected $\DGA$ is said to be \emph{formal} if its minimal model is weakly equivalent to the minimal model of its cohomology algebra. Recall that a {\it minimal algebra} $m$ is a free graded algebra which can be written as an increasing union of $m_i$'s where $m_0=m_1=\QQ$, and $m_{k-1} \to m_k$ is a Hirsch extension for every $k$ (\cite[Theorem 10.3]{PGJM}). A {\it minimal model} of a $1$-connected $\DGA$ $u$, is a pair $(m,\rho)$, where $m$ is a minimal algebra and $\rho:m\to u$ is a quasi-isomorphism.

A simply connected space is said to be \emph{formal} if the corresponding minimal model is formal. The rational homotopy groups of formal spaces can be computed from its cohomology algebra and by rational Postnikov tower (\cite[Theorem 3.3]{DGMS}). A graded algebra $A^{\ast}$, is said to be $k$-intrinsically formal, if for any minimal algebra $m$ with $H^{\ast}(m)\cong A^{\ast}$ the sub $\DGA$ $m_k$ of $m$ generated by elements degree $\leq k$, is unique up to isomorphism. If $m_k$ is unique for every $k$, the graded algebra $A^{\ast}$ becomes intrinsically formal, and any space with cohomology algebra $A^{\ast}$ will be formal. In this section, we introduce the notion of formality and intrinsic formality in the equivariant setup with some examples.  We first prove the following facts in the non-equivariant case.  

\begin{prop}\label{prop:kiso implies quasi iso}
    Let $m,m'$ be two minimal algebras with $m_k\cong m'_k$. If there is a morphism $m_{k+1}\to m'_{k+1}$ which is a $(k+1)$-isomorphism, then $m_{k+1}\cong m'_{k+1}$. 
\end{prop}

\begin{proof}

    Let $f:m_{k+1}\to m'_{k+1}$ be a  $(k+1)$-isomorphism, which is an extension of the isomorphism from $m_k\to m_k'$. 

    Consider the diagram 
 \[
\xymatrix{
m_k\cong m'_k \ar[d]\ar[r]& m_{k+1}\ar[d]^{f} \\
m'_{k+1}\ar@{-->}[ru]^g\ar[r]^{id} & m'_{k+1} } 
\]

where the bottom horizontal arrow is the identity. The obstructions to finding a lift $m_{k+1}' \to m_{k+1}$ successively lie in the relative cohomology (\cite[Proposition 11.1]{PGJM})  $H^{i+1}(m_{k+1}, m'_{k+1}, V^i)$  where $V^i$ are the degree $i$ generators of $m'_{k+1}$.  If the relative cohomology $H^{i+1}(m_{k+1}, m'_{k+1})$ vanishes for $i+1 \leq  k+2$, then all the obstructions vanish.  Note, we only have to consider $i$ up to $k+2$, as $m_{k+1}$ is generated by elements of  degrees $ \leq k+1$.

 Consider the long exact sequence in cohomology, 

$$\cdots \to H^k(m_{k+1}) \to H^k(m_{k+1}') \to  H^{k+1}(m_{k+1}, m_{k+1}') \to  H^{k+1}(m_{k+1}) \to \cdots.
$$

Now, from our assumptions on the map $f$, it follows that $H^{\leq  k+2}(m_{k+1}, m_{k+1}')$ vanishes. Therefore, we have a lift $g$ such that $fg$ is homotopic to the identity on $m_{k+1}'$. 

The existence of this $g$ implies that $f$ is cohomologically surjective in all degrees. Also, $g$ is cohomologically injective in all degrees.  Since $f$ is a cohomological isomorphism in all degrees $\leq k+1$, so is $g$.

Similarly, consider  the diagram 

 \[
\xymatrix{
m_k\cong m_k' \ar[d]\ar[r]& m'_{k+1}\ar[d]^{g} \\
m_{k+1}\ar[r]^{id}\ar@{-->}[ru]^{h} & m_{k+1}. } 
\]
We  get a lift $h$ as before, and conclude that $gh$ is homotopic to the identity. Thus, $g$ is cohomologically surjective in all degrees. Therefore, it is a cohomological isomorphism in all degrees.

Since $fg \simeq id$, $f$ is also a  cohomological isomorphism. 

Both $m_{k+1}$ and $m'_{k+1}$ are  minimal algebras, and  quasi-isomorphism implies isomorphism. Hence $m_{k+1}\cong m'_{k+1}$.

\end{proof}

\begin{corollary}
   Let $m$ and $m'$ be  minimal algebras. A $k$-isomorphism between $m_k$ and $m'_k$ induces a quasi-isomorphism.

\end{corollary}

We can talk about formality and intrinsic formality in the equivariant case as follows. 

\begin{definition}\label{def:formal}
 We say a system of $\DGA$s, $(\sA,d)$ is equivariantly formal if there is a weak equivalence between $(\sA,d)$ and the injective envelope of $(H(\sA), 0)$. That is, the isomorphism of graded algebras is realized by a zig-zag of quasi-isomorphisms of the system of $\DGA$s.  A $G$-space $X$ is said to be equivariantly formal if the minimal system of $\DGA$s representing it is equivariantly formal.  

 \end{definition}

\begin{definition}[Equivariantly intrinsically formal]\label{def:I.F}
A diagram of graded algebras $\cA$ over $\sO_G$, for some finite group $G$ is called equivariantly $k$-intrinsically formal (abbreviated to equivariantly $k$-I.F.) if for any two minimal systems $\sM$ and $\sM'$ with $H^\ast(\sM)=\cA=H^{\ast}(\sM')$, the sub-systems $\sM_k$ and $\sM'_k$ have the property that, there is a map from $\sM_k \to \sM'_k$ or $\sM_k'\to \sM_k$ which is a $k$-isomorphism.
 For any two such minimal systems $\sM$ and $\sM'$ if there is a map $\phi:\sM \to \sM'$ which is a  quasi-isomorphism, we say $\cA$ is equivariantly intrinsically formal (abbreviated to equivariantly I.F.).
\end{definition}

\begin{remark}
\begin{enumerate}
\item The above notion of intrinsic formality is compatible with that of the non-equivariant case by Proposition \ref{prop:kiso implies quasi iso}.

    \item 
Any diagram of graded algebras is $2$-intrinsically formal.
Let $\cA$ be a diagram of graded algebras and let $\sM$ be any minimal algebra such that $H^\ast(\sM)=\cA$. By the definition of minimal algebra $\sM_{0}=\mathbb{Q}=\sM_{1}$. For any diagram of vector spaces $\underline{V}$, $H^3(\sM_{1};\underline{V})=0$. So in order to get $\sM_2=\sM_{1}^{\gamma}(\underbar{V})$, the choice of $\gamma\in H^3(\sM_{1};\underbar{V})$ is unique and $\sM_2$ is uniquely determined.

Thus, if $\cA$ is  equivariantly I.F. and the diagram of cohomology algebras for a $G$-space $X$ is isomorphic to   $\cA$, then $X$ is equivariantly formal. 

\end{enumerate}
\end{remark}

\begin{prop}\label{prop:IIF}
    Let a cohomology diagram $\cA$ be equivariantly $n$-I.F. with $\cA^{n+1}=\cA^{n+2}=0$. Let $\sM$ be the minimal model for $\cA$ with $\sM_n=\sM_{n+1}$. Then $\cA$ is equivariantly $(n+1)$-I.F.  
\end{prop}
\begin{proof}
    First, note that if $\sM$ is any minimal system then for any $n\geq 0$ the inclusion $\sM_n \to \sM$ is an $n$-isomorphism. 

    Let $\sM'$ be another minimal system with $H^{\ast}(\sM')=\cA$. As $\cA$ is equivariantly $n$-I.F., let  $\phi_n:\sM_n \to \sM'_n$ be an $n$-isomorphism. 

Since $\sM_n=\sM_{n+1}$, we rewrite the map as  $\phi_n:\sM_{n+1} \to \sM'_{n+1}$.

Note that $H^{n+2}(\sM_{n+1})=H^{n+2}(\sM_n)=0$ as $\sM_{n+1} \to \sM$ is $(n+1)$-isomorphism and $H^{n+1}(\sM)=\cA^{n+1}=0$. We claim that $H^{n+2}(\sM'_{n+1})=0$. Now if $\alpha \in H^{n+2}(\sM_{n+1}')$ is non-zero then as $\sM'_{n+1}\to \sM'$ is an $(n+1)$ isomorphism by injectivity, we have a non-zero member which is the image of $\alpha$ in $H^{n+2}(\sM')\cong \cA^{n+2}=0$ and this is a contradiction. The map $H^{n+2}(\sM_{n+1})\to H^{n+2}(\sM'_{n+1})$ induced by $\phi_n$ is a zero map and hence a monomorphism. 

Using a similar argument, we conclude that $H^j(\sM_{n+1}) \to H^j(\sM'_{n+1})$ is an $n$-isomorphism for $j\leq n$.

Hence, $\cA$ is equivariantly $(n+1)$-I.F.

\end{proof}

\begin{exmp}
\begin{enumerate}
    \item Consider the space $X=S^n\vee S^n \vee \cdots \vee S^n$, which is  $(p+1)$ many copies of spheres of dimension $n\geq 2$. There is an action of $C_p$ on $X$ which permutes the first $p$ copies of $S^n$ and keeps the last copy fixed. Note that the cohomology algebra is equivariantly $n$-I.F. as the fixed point sets are $(n-1)$ connected. The cohomology algebra also satisfies the hypothesis of Proposition \ref{prop:IIF} so it is equivariantly $(n+1)$-I.F. Since all higher cohomologies are zero, by Proposition \ref{prop:IIF} we conclude that the cohomology algebra is equivariantly I.F.  

    \item Consider $X=S^n \times \cdots \times S^n$ be product of $p$-copies of $S^n$'s with $p$ a prime number and $n\geq 3$. Then there is a $C_p$-action on $X$ by permutation and the fixed point set is homeomorphic to $S^n$. As both fixed point sets are $(n-1)$ connected, so the space is equivariantly $n$-I.F. As $n\geq 3$, using Proposition \ref{prop:IIF} the cohomology diagram is equivariantly $(n+1)$-I.F. 
\end{enumerate}
  
\end{exmp}

We give a condition under which a map of systems of $\DGA$s can be extended to a map from an elementary extension of the domain.

\begin{prop}\label{prop:extension}
 Let  $f:\mathcal{U}\to \mathcal{B}$ be a map of systems of $\DGA$s and $\mathcal{U}^{\alpha}(\underline{V})$ is an elementary extension with respect to some $\alpha$. 
 \begin{enumerate}
     \item 
 
 If $f':\underline{V}\to \mathcal{B}$ satisfies $f\alpha=df'$, one can extend $f$ to $\mathcal{U}\otimes \mathbb{Q}(\underline{V})$ using $f'$ on $\underline{V}$.
   \item  The converse is also true. Let $f:\mathcal{U}^{\alpha}(\underline{V})\to \mathcal{B}$ is a morphism, i.e., $df=fd$. If $f'=f|_{\underline{V}}$, then $f'$ satisfies $f\alpha=df'$.
 
 \end{enumerate}
\end{prop}

\begin{proof}

Suppose $f':\underbar{V}\to \cB$ satisfies $f\alpha=df'$, we can extend $f$ to $\cU\otimes \QQ(\underbar{V})$ using $f'$ on $\underbar{V}$; the condition on $f'$ ensures that this map respects the differential. We now extend the map to $\underbar{V}_0$ using the injectivity of $\cB$ and define it on the rest of the resolution inductively. Given $f'$ on $\underbar{V}_i$, we must define a map $f':\underbar{V}_{i+1}\to \cB$ such that $f'd=df'$. To ensure this is satisfied, we consider the differential from $\underbar{V}_i$ to $\underbar{V}_{i+1}$; defined by $d=(-1)^i\alpha_i+w_i$. We  need to find a  map $f'$ such that $$f'((-1)^i\alpha_i+w+i)=(-1)^if'\alpha_i+f'w_i=df',$$

or equivalently, $(-1)^if'\alpha_i-df'=fw_i$. Observe that since 
\begin{align*}
((-1)^if'\alpha_i-df)w_{i-1} &=(-1)^if\alpha_iw_{i-1}-dfw_{i-1}\\
&=(-1)^ifd\alpha_{i-1}-fdw_{i-1}\\
&=fd((-1)^i\alpha_{i-1}-w_{i-1})\\
&=fd(-d)\\
&=0
\end{align*}

the map $(-1)^if'\alpha-df$ vanishes on $\Im(w_{i-1})\subset \underbar{V}_i$, and we have

 \[
\xymatrix{
\underline{V}_i/\Im w_{i-1} \ar@{^{(}->}[r] \ar[d] & \underline{V}_{i+1}\ar@{-->}[dl]^{f'}\\
{\cB}^{n+i+1}\\
} 
\]

where, the map $ \underbar{V}_i/\Im w_{i-1}\to \cB^{n+i+1}$ is given by $(-1)^if\alpha_i-df$. Since $\cB$ is injective,  we can define $f'$ on $V_{i+1}$. Continuing in this manner, we extend $f'$ to all generators and therefore, to a $\DGA$ map on all of $\cU(\underbar{V})$.

Conversely, take $f'=f|_{\underline{V}}.$ Consider the injective resolution of $\underline{V}$, i.e., $$\underline{V}\to \underline{V_0}\to \underline{V_1}\to \cdots.$$ The map $\underline{V_0}\xrightarrow{w_0} \underline{V_1}$ when restricted to $\underline{V}(\subset \underline{V_0})$ is $w_0|_{\underline{V}}=0,$ since  $\underline{V_1}$ is the injective envelope of  $\coker (w)$. By definition, if we pick an element $x\in \underline{V}$, the derivation $d$ has no horizontal component on $\underline{V}$, that is,  $d=\alpha$ on $\underline{V}$. Hence $f\alpha=df'$, on $\underline{V}$,  where $f'=f|_{\underline{V}}$. 

\[
\xymatrix{
\underline{V}\ar[d]_{\alpha} \ar[r]^{w} & \underline{V}_0\ar[d]^{ \alpha_0}\ar[r]^{w_0}& \underline{V}_1\ar[d]^{\alpha_1}\ar[r]^{w_1} & \underline{V}_2\ar[r]\ar[d]^{\alpha_2} & \cdots  \\
\underline{Z}^{n+1}(\mathcal{U})\ar[r]         & \mathcal{U}^{n+1}\ar[r]_{d}&\mathcal{U}^{n+2}\ar[r]_d & \mathcal{U}^{n+3}_d \ar[r] &\cdots
}
\] 
\end{proof}

\indent We now recall the construction of the minimal model of a system of $\DGA$s from \cite[Thm. 3.11]{LS02}. Let $\cU$ be a system of $\DGA$s then inductively we build the minimal model $\sM=\cup\sM_d$, where each $\sM_n$ is an  elementary extension as in Section \ref{section:background}, so that $\sM_n=\sM_{n-1}^{\gamma}(\underbar{V})$, here $\gamma \in H^{n+1}(\sM_{n-1};\underbar{V})$. 

Consider the diagram:

 \[
\xymatrix{
        & \sM_{n-1}\ar[rd]^{\rho} \ar[d]_{\alpha} \\
\ker{\beta} \ar[r]_{i}       & \sM'_{n-1}\ar[r]_{\beta} & \cU } 
\]

Here $\rho$ is an $(n-1)$-isomorphism, $\alpha$ is a quasi-isomorphism and $\beta$ is surjective.

The system of $\DGA$s, $\sM'_{n-1}=\sM_{n-1}\otimes (\otimes_H\QQ(\underbar{U}^\ast_H\oplus \sum \underbar{U}^\ast_{H}))$, where $U_H,\underline{U}_H^{\ast}$ are from \Cref{equation:24} and $\sum \underbar{U}^\ast_H$ is the system obtained from $\underline{U}_H^{\ast}$ by considering the  degree of every element shifted by $+1$. The differential on $\sM'_{n-1}$ is defined accordingly. Here, the map

$$\beta|_{\sM_{n-1}}=\rho, ~ \beta|_{\underbar{U}}=\text{id} \text{ and } \beta(\sum x)=dx, \text{ for } x\in \sum \underbar{U}.$$

Let $\underline{\QQ}$ be the coefficient system defined by $\underline{\QQ}(G/H)=\QQ$ for every subgroup $H$ of $G$. Define $R:=\ker\beta \oplus \underline{\QQ}$, and the diagram of vector spaces $\underbar{V}=H^{n+1}(R)$. 
The map $\gamma$ is the elementary extension  obtained by considering $[id] \in H^\ast(\underbar{V};\underbar{V})$ and considering its image, under the inclusion $i:\ker\beta \to \sM'_{n-1}$, in $H^\ast(\sM'_{n-1};\underbar{V})$. Since $\alpha$ is quasi-isomorphism, there is a pre-image $\gamma\in H^\ast(\sM_{n-1};\underbar{V})$, of $i^\ast[id]$ such that $\alpha^\ast[\gamma]=i^\ast[id]$.

Let $X$ be a $G$-space and  $\cA$ be its cohomology diagram. Let $\sM$ be the minimal model of $\cA$.  Then there exists a quasi-isomorphism $\sM\xrightarrow{\rho} \cA$.  The system of $\DGA$s, $\sM$, is a minimal system and by definition we have $\sM=\cup_{i\geq 0}\sM_i$, where $\sM_n=\sM_{n-1}^{\gamma}(\underbar{V})$. 

\begin{definition}\label{def:associ:vect}
Let $\cA$ be an injective cohomology diagram and $(\sM,\rho)$ be its minimal model. Let $\sM_n=\sM_{n-1}^{\gamma}(\underbar{V})$, be the $n$-th stage construction of $\sM$, which is obtained by taking elementary extension of $\sM_{n-1}$ with the injective resolution of $\underbar{V}$. We refer to $\underbar{V}$ as the  \emph{$n$-th stage associated diagram of vector spaces of $\cA$}. 
\end{definition}

Let $\sU$ be a system of $\DGA$s over $\sO_G$ and  $\sM$ be the minimal model of $\sU$. Let the restriction at the $(n-1)$th level, $\rho: \sM_{n-1}=\sN \to \sU$ be such that $\rho $ is an $(n-1)$-isomorphism, $\rho^\ast[\gamma]=0$, for $\gamma \in H^{n+1}(\sN; \underline{V})$ and $\sN^{\gamma}(\underline{V})$ is the $n$-th stage of $\sM$. That is, $\rho: \sM_{n}=\sN^{\gamma}(\underline{V})\to \sU$ is an \textbf{$n$-isomorphism}, a cohomology isomorphism up to degree $n$ and monomorphism at degree $n+1$.\\

Let $\gamma'\in H^{n+1}(\sN;\underline{V})$. 
We say $\sN^{\gamma'}(\underline{V})$ satisfies \textbf{Condition $\C_n$} with respect to $\sN^{\gamma}(\underline{V})$ if the following is true;  
\begin{enumerate}

    \item[] With the above assumptions, there exists $\kappa \in Aut(\sN)$ and $g\in Aut(\underline{V})$ such that $\kappa {\gamma}=d_{\gamma'}(g+\beta)$, where $\beta :\underline{V} \to \sN$.
\end{enumerate}

\begin{remark}

Note that if we choose $\sigma=\rho^\ast$ and $\gamma=\gamma'$ then $\rho^\ast[\gamma]=0$. Also, there always exists $\kappa=id\in Aut(\sN)$ and $g=id\in Aut(\underline{V})$ so that $\kappa\circ id=id\circ \kappa $ (considering $\beta=0$).

\end{remark}

\begin{prop}\label{prop:plural}

Let $\sM$ be a minimal system with $\sM_{n-1}=\sN$ and $\sM_n=\sN^{\gamma}(\underline{V})$. Assume $\gamma'\in H^{n+1}(\sN; \underline{V})$.
Then the following statements are equivalent. 
\begin{enumerate}

    \item $\sN^{\gamma'}(\underline{V})$ satisfies \textbf{Condition} $\mathbf{C_{n}}$. 
    \item  $\sN^{\gamma}(\underline{V})$ and $\sN^{\gamma '}(\underline{V})$ are isomorphic.
\end{enumerate}
\end{prop}

\begin{proof}

$(2) \to (1):$  If $f: \sN^{\gamma}(\underline{V})\to \sN^{\gamma'}(\underline{V})$ is an isomorphism then by Remark \ref{rmk:isolevelwise} and Lemma \ref{lemm:13.57}, we may assume that $f$ takes $(\sN^{\gamma}(\underline{V}))(n)$ to $(\sN^{\gamma'}(\underline{V}))(n)$ that is, $f$ is a level-wise isomorphism.   
Since isomorphisms are quasi-isomorphisms, and a quasi-isomorphism between minimal systems is homotopic to an isomorphism which is level invariant by Remark \ref{rmk:isolevelwise}, we can assume, $f:\sN \to \sN$ is an isomorphism.  For $x\in \underline{V}$ we have $f(x)=x+b(x)$, where $b:\underline{V}\to \sN$.

As, $fd=df$,
 we have for $x\in \underline{V}$, $f{\gamma} (x)=d_{\gamma'} f(x)=d_{\gamma'}(x+b)=d_{\gamma'}(id(x)+b)$ (here the first equality comes from the converse part of Proposition \ref{prop:extension}). So the maps $f:\sN \to \sN$ and $id:\underline{V}\to \underline{V}$ give that Condition $\C_n$ is satisfied. 

$(1) \to (2)$ If there exists $\kappa \in Aut(\sN)$ and $g \in Aut(\underline{V})$ such that $\kappa {\gamma}=d_{\gamma'}(g+\beta)$ then by Proposition \ref{prop:extension} this implies that, $\kappa$ extends to a map $\Tilde{\kappa}:\sN^{\gamma}(\underline{V})\to \sN^{\gamma'}(\underline{V})$. By Proposition \ref{prop:11.53}, $\Tilde{\kappa}$ is isomorphism.

\end{proof}

\begin{exmp}

We define a $G=C_2$ action  on $X=S^3\vee S^3 \vee S^5$, where $C_2$ acts on $(S^3 \vee S^3)$ by switching copies and acts on $S^5$ trivially.

The fixed point set $X^G=S^5$. We denote the cohomology diagram of $X$ by $\cA$ and, the generators of the cohomology by $x$, $y$, and $z$. Note that, $deg(x)=3=deg(y)$ and $deg(z)=5$. 

Then $\cA(e)=\mathbb{Q}(x,y,z)/D$ and $\cA(G)=\mathbb{Q}(z)/E$, where $D=<x^2, z^2, xy,xz,yz>$ and $E=<z^2>$.

We will compute the minimal system for $X$ under $G$ action up to the $6$th stage.
We get that, 

$$\sM_0=\mathbb{Q}=\sM_1=\sM_2,$$

\[
\sM_3=
\begin{cases}
     \wedge(a_3,b_3), & \text{at G/e}\\
    \mathbb{Q},  & \text{ at G/G}
\end{cases}
\] 
 $$\sM_4 \cong \sM_3,$$

\[
\sM_5=
\begin{cases}
     \wedge(a_3,b_3,c_5,d_5), & \text{at G/e}\\
    \wedge(c_5),  & \text{ at G/G}
\end{cases}
\] 

The differential is $d(a_3)=0=d(b_3),$ $d(d_5)=ab$ and $d(c_5)=0$, so

$$\sM_6\cong \sM_5.$$

At each stage, $H^{n+2}(\sM_n;\underbar{V})$ is zero for $n\leq 5$, thus, Condition $C_5$ is satisfied trivially .

Later in Example \ref{example sec5}(\ref{wedgeeeodd}), we show that this diagram of graded algebras is equivariantly $5$-I.F.

\end{exmp}

\begin{exmp}\label{example:plural}
Consider $X=(S^3\vee S^3)\times S^5$ with a diagonal action of $G=C_2$, where $C_2$ acts on $(S^3 \vee S^3)$ by switching copies and acts on $S^5$ trivially. Then $X^G\cong S^5$. 
We denote the cohomology diagram by $\cA$, the generators of the cohomology by $x$, $y$ and $z$. Note that, $deg(x)=3=deg(y)$ and $deg(z)=5$, so 

\[
\cA:=
\begin{cases}
    \wedge(x,y)/<xy>\otimes \wedge (z), & \text{at G/e}\\
    \wedge (z),  & \text{ at G/G}.
\end{cases}
\]

At the $7$th stage, we have, 

\[
\sM_7= \sM_6^{\gamma_1}(\underbar{V})
\begin{cases}
     \wedge(a_3,b_3,c_5,d_5,e_7,f_7), & \text{at G/e}\\
    \wedge(c_5),  & \text{ at G/G}
\end{cases}
\] 
with $d(a)=d(b)=d(c)=0,$, $d(d_5)=ab$, $d(e_7)=ad$, $d(f_7)=bd$, 
and
\[
\sM'_7= \sM_6^{\gamma_2}(\underbar{V})
\begin{cases}
     \wedge(a_3,b_3,c_5,d_5,e'_7,f'_7), & \text{at G/e}\\
    \wedge(c_5),  & \text{ at G/G}
\end{cases}
\] 

with $d(a)=d(b)=d(c)=0,$, $d(d_5)=ab$, $d(e'_7)=ad+ac$, $d(f'_7)=bd$. Clearly, $\sM_7$ and $\sM'_7$ are not quasi-isomorphic as a system of $\DGA$s by Proposition \ref{prop:plural}.
\end{exmp}

Let $\gamma, \gamma', \underline{V}$ be as in \Cref{prop:plural}.  We say $\sN^{\gamma'}(\underline{V})$ satisfies \textbf{Condition $(\PH)_n$} if the following holds. 
\begin{enumerate}

\item $\sN^{\gamma'}(\underline{V})$ does not satisfy Condition $\C_n$ with respect to $\sN^{\gamma}(\underline{V})$. 
    \item 

There exists a map $\sigma_n^{\ast}: H^{\ast}(\sN^{\gamma'}(\underline{V}))\to \mathcal{U} $, such that $\sigma_n^i$ is isomorphism for $i\leq n$ and monomorphism for $i=n+1$ and 

\[
 \xymatrix{
H^{\ast}(\sN) \ar[d]_{i^{\ast}}\ar[r]^{\rho^{\ast}} &\mathcal{U}\\
H^{\ast}(\sN^{\gamma'}(\underline{V}))\ar[ur]^{\sigma_n^{\ast}}
}
\]
commutes for $\ast \leq n-1$.

\end{enumerate}

Inductively, if the Condition $(\PH)_{n}$ holds for all $n\geq 1$, we say that the \emph{Condition plural homotopy type} is satisfied. In this case we will get a new minimal system $\sM_{\infty}=\sN\cup\sN^{\gamma'}(\underline{V})\cup (\sN^{\gamma'}(\underline{V}))_{\delta}(\underline{W})\cup \cdots $, and a map $\sigma_{\infty}: H^{\ast}(\sM_{\infty}) \to \mathcal{U}$ which is an isomorphism. Thus, we get more than one non-isomorphic minimal algebra with the same cohomology algebra.

\begin{remark}
\begin{enumerate}

\item It is difficult to check \textbf{Condition} $\mathbf{C_n}$ for each $n$, as it involves the computation of elementary extension at each level.

\item  Note, this condition allows us to find plural homotopy types. But when this condition is not satisfied, it is not clear  whether the given diagram of graded algebras will be equivariantly intrinsically formal or not. We now provide a few examples where the Condition plural homotopy type is not satisfied. We prove later in the article that they are equivariantly intrinsically formal.

\end{enumerate}
\end{remark}

We give a couple of examples of  equivariantly intrinsically formal and equivariantly $n$-intrinsically formal diagrams of graded cohomology algebras.

\begin{exmp}\label{odd:product}
Let $n $ be an odd integer $\geq 3$ and $p$ be a prime.
Let $X= S^n \times \cdots \times S^n$ be the product of $p$ copies of $S^n$ with the  $G=C_p$ action $t(a_1,a_2,\cdots, a_p)= (a_2,a_3,\cdots, a_{p-1},a_1)$ where $t$ is a generator of $C_p$. The fixed points of $X$ under $G$, $X^G \cong S^n$. Thus, the cohomology diagram of $X$, which we denote by $\cA$, is given by 
$$ \cA =
\begin{cases}
\QQ_3[x_1,x_2,\cdots, x_p]/<x_i^2~|~ i=1,2,\cdots, p> & \text{ at G/e} \\
\QQ_3[y]/<y^3>  & \text{ at G/G},
\end{cases}$$

where $y$ corresponds to the generator of the cohomology algebra for $X^G$. We want to compute the minimal model $\sM$ for this cohomology diagram.

On further calculation we get, $$\sM_0=\cdots =\sM_{n-1}=\underline{\QQ}$$ and 
\[
\sM_n=
\begin{cases}
     \wedge(a_{12},a_{23},\cdots a_{(p-1)p},b), & \text{at G/e}\\
    \wedge(b),  & \text{ at G/G}.
\end{cases}
\]

with zero differential. Here the process ends at the $n$-th stage, since we get a quasi-isomorphism from $\sM_n \to \cA$. So the minimal system is obtained at the $n$-th stage.

From the calculation, we see that $$H^{r+1}(\sM_{r-1})=0$$ for every $r$, and we conclude that Condition $(PH)_r$ is not satisfied for any $r$. Also, if $\sM'$ is any minimal algebra with $\underbar{H}^{\ast}({\sM'})=\cA$, then $\sM'$ consists of at least as many generators as $\sM_n$. Thus, one can define an $n$-isomorphism via inclusion from $\sM_n\to \sM'$. Later in  \Cref{example sec5}(\ref{oddddd}) we will show that the cohomology diagram is equivariantly I.F. 
\end{exmp}

\begin{exmp}\label{even:product}
Let $n $ be an even integer $\geq 2$ and $p$ be a prime.
Let $X= S^n \times \cdots \times S^n$ be the product of $p$ copies of $S^n$ with   $G=C_p$ action given by $t(a_1,a_2,\cdots, a_p)= (a_2,a_3,\cdots, a_{p-1},a_1)$ where $t$ is a  generator of $C_p$. The fixed points of $X$ under $G$ is $X^G \cong S^n$. Thus, the cohomology diagram of $X$, which we denote by $\cA$, is given by 
$$ \cA =
\begin{cases}
\QQ_n[x_1,x_2,\cdots, x_p]/<x_i^2~|~ i=1,2,\cdots, p> & \text{ at G/e} \\
\QQ_n[y]/<y^2>  & \text{ at G/G},
\end{cases}$$

where $y$ corresponds to the generator of the cohomology algebra for $X^G$. We want to compute the minimal model for the cohomology diagram. 

Putting all the stages together, we get the following:

 $$\sM_0=\underline{\mathbb{Q}}=\sM_1=\cdots =\sM_{n-1}.$$

\[
\sM_n=
\begin{cases}
     \wedge(a_{12},a_{23},\cdots a_{(p-1)p},b), & \text{at G/e}\\
    \wedge(b),  & \text{ at G/G}.
\end{cases}
\] 

Since the both $G/e$ and $G/G$ levels of the spaces are $(n-1)$-connected, the cohomology diagram is equivariantly $(n-1)$-I.F. 

Next, we claim that the cohomology diagram is equivariantly $n$-I.F. If we consider any minimal system $\sM'$ with $H^{\ast}(\sM')=\cA$ then $\sM'_{n-1}=\sM_{n-1}$. Now $\sM'_{n}=\sM_{n-1}(\underbar{V})$. Then we claim that there is a map $\sM_n \to \sM'_n$ that is an $n$-isomorphism. 

First, note that $\sM(G/e)$ and $\sM(G/G)$ are (non-equivariantly)  minimal algebras. For any other minimal system $\sM'$ and, for any $H\leq G$,  $\sM'_n(G/H)$ has at least as many generators as that in  $\sM_n(G/H)$. Since the cohomology diagrams of both $\sM$ and $\sM'$ are isomorphic, the generators which contribute to the non-zero cohomology classes of $\sM_n$ mapping to the generators of $\sM'_n$ which are non-trivial classes defines an inclusion map $\sM_n\to \sM'_n$. This map is an $n$-isomorphism by construction.  

Also,
$$\sM_n=\cdots =\sM_{2n-2}.$$

The map $\sM_{2n-2}=\sM_n\to \sM'_n\to \sM'_{2n-2}$ is a $(2n-2)$-isomorphism since $\cA^i=0$, for $n+1\leq i\leq 2n-1$.

\[
\sM_{2n-1}=
\begin{cases}
     \wedge(a_{12},a_{23},\cdots a_{(p-1)p},c_{12}, \cdots ,c_{(p-1)p} ,b,b'), & \text{at G/e}\\
    \wedge(b,b'),  & \text{ at G/G}.
\end{cases}
\] 

with $d(c_{ij})=a_{ij}^2$ for all $i,j$ and $d(b')=b^2$.

Given that  $\sM_{2n-1}(G/H)$ are minimal algebras for all $H<G$, using the earlier argument (the way we show that $\cA$ is equivariantly $n$-I.F.) we can show that  $\cA$ is equivariantly $(2n-1)$-I.F.

Later in Example \ref{example sec5}(\ref{evennnn}) we will show that the product of even spheres under the above action is equivariantly formal.  

\end{exmp}

\section{Injectivity of the associated diagram of vector spaces}

In the non-equivariant case, various authors (\cite[Theorem 3.2]{SY82}, \cite[Example 6.5]{HARST}, \cite{HIRO}) give different methods to  determine the plural homotopy types of a given graded algebra.
These results often use the fact that at the $n$-th extension stage, we are adding generators only in degree $n$.  This  is not true in the equivariant case. However, if the associated diagram of graded vector spaces at the $n$-th stage  is  injective then, the generators added are only in degree $n$.  

In this section, we restrict our group to $G=C_p$, for prime $p$, and describe conditions under which a diagram of $\DGA$s or a diagram of vector spaces is injective.

\begin{prop}\label{prop:injcohom}
Let $G=C_p$, where $p$ is a prime. Let $\sA\in \DGA^{\sO_G}$. Then  $\sA$ as an element of $Vec^{\ast}_G$ is injective if and only if the map $\sA(\hat{e}_{e,G}): \sA(G/e)\to \sA(G/G)$ is surjective. 
\end{prop}
\begin{proof}
The injective envelope for $\sA$ is given by $\sI(\sA)=\cI_G^\ast \oplus \cI_e^\ast$, where $\cI_e^{\ast}$ and $\cI_G^{\ast}$ are systems corresponding  to the vector spaces $I_e=\ker\sA(\hat{e}_{e,G}) $ and $I_G=\sA(G/G)$ respectively. Given that the map $\sA(\hat{e}_{e,G})$ is surjective, we have a short exact sequence 
\beq\label{ses1}
0\to \ker\sA(\hat{e}_{e,G}) \to \sA(G/e) \to \sA(G/G)\to 0,
\eeq
which splits as $\QQ$-vector spaces. So $\sA(G/e)=\ker\sA(\hat{e}_{e,G})\oplus \sA(G/G)$.

Also,  
\begin{eqnarray}
\cI^\ast_G(G/G ) &= & Hom_{\mathbb{Q}(G/G)}(\mathbb{Q}(G/G)^G,I_G)\cong I_G,\\
\cI^\ast_G(G/e)&=& Hom_{\mathbb{Q}(e)}(\mathbb{Q}(G/G)^e, I_G)\cong I_G,\\
\cI^\ast_e(G/G)&=&Hom_{\mathbb{Q}(G)}(\mathbb{Q}(G/e)^G,I_e)=0,\\
\cI^\ast_e(G/e)&=& Hom_{\mathbb{Q}(G/e)}(\mathbb{Q}(G/e)^e, I_e)\cong I_e. \label{ses2}
\end{eqnarray} 

Since, $\sI(A)= \cI_G^\ast \oplus \cI_e^\ast$ we get that  
$$\sI(\sA)(G/e)=\cI^\ast_G(G/e)\oplus \cI^\ast_e(G/e)= \ker\sA(\hat{e}_{e,G})\oplus \sA(G/G)= \sA(G/e),$$
 and 
 $$\sI(\sA)(G/G)=\cI^\ast_G(G/G)\oplus \cI^\ast_e(G/G)= I_G=\sA(G/G). $$

 Therefore, the injective envelope of $\sA$ is itself implying that  $\sA$ is injective.
 
 Conversely, if $\sA$ is injective, then $\sA\cong \cI_G^\ast \oplus \cI_e^\ast$. Then using \ref{ses1} $-$ \ref{ses2} we get that  the map $\sA(\hat{e}_{e,G})$ is the projection $$\sA(G/G)\oplus ker\sA(\hat{e}_{e,G})= \cI^\ast_G(G/e)\oplus \cI^\ast_e(G/e) \to \cI^\ast_G(G/G)\oplus \cI^\ast_e(G/G)=\sA(G/G)$$ and hence is surjective.  
\end{proof}

\begin{exmp}\label{exmp:noninjective}
Consider the $G$-space $X=S^3$, where $G=C_2$ acts on $S^3$ by reflection, which fixes the equator sphere $S^2$. So here $G=C_2$, $X^G=S^2$ and $X^e=S^3$. The corresponding cohomology diagram is given by $H^\ast(X;\underline{\QQ})$, which is not injective. This follows from Proposition \ref{prop:injcohom}.  
\end{exmp}

Note that if a cohomology diagram $\cA$ is injective and for each $n$, and the associated diagram of vector spaces for $\cA$ is injective, then at elementary extension we only add elements of degree $n$ to  $\sM_{n-1}(G/H)$ to obtain $\sM_n(G/H)$, for all $H\leq G$. Then by \cite[Lemma 3.2]{DGMS} the differential will be level-wise decomposable.  Also, the map $\rho(G/H): (\sM(G/H),d) \to (\cA(G/H),0)$ is a quasi-isomorphism and surjective, since the differential in $\cA(G/H)$ is $0$. Thus, by the Lifting Lemma \cite[Lemma 12.4]{FHT}, there is a map $\alpha: N_H \to \sM(G/H)$. Since $\alpha_H$ is a quasi-isomorphism between two minimal algebras, it is an isomorphism. In view of this, we have the following Proposition.

\begin{prop}\label{prop: isonon-eq and equivariant}
Let $\cA$ be a cohomology diagram over $\sO_G$, which is injective, and for each $n$, the associated diagram of vector spaces of $\cA$ is injective. If $(N_H, \rho_H)$ is the minimal model for $\cA(G/H)$ where  $\rho_H: N_H\to \cA(G/H)$ is a quasi-isomorphism and $\sM$ be a minimal system for $\cA$, then $\sM(G/H)\cong N_H$.  
\end{prop}

\begin{remark}
Let us consider the above proposition when  $G=C_p$ for $p$ prime and $(N_e,\rho_e)$ and $(N_G,\rho_G)$ are minimal models for $\cA(G/e)$ and, $\cA(G/G)$ respectively. If there exists a  map $\theta: N_e \to N_G$ such that the following diagram commutes   
\[
\xymatrix{
N_e\ar[r]^{\alpha_e}\ar[d]_{\theta^\ast}  & \sM(G/e) \ar[d]^{\sM(\hat{e}_{e,G})}\ar[r]^{\rho_e} & \cA(G/e)\ar[d]^{\cA(\hat{e}_{e,G})}\\
N_G \ar[r]_{\alpha_G}   & \sM(G/G) \ar[r]_{\rho_G}    & \cA(G/G)} 
\]

Then the minimal system for the diagram of graded algebras $\cA$ can be given by $$\sM(G/e)=N_e,~~\sM(G/G)=N_G.$$
\end{remark}

 We have the following result. 
\begin{prop}\label{Z_P injimpliesminimal}
Let $G=C_p$, for $p$ prime. If the structure map in the cohomology diagram $\cA$  $\cA_{e,G}: \cA(G/e)\to \cA(G/G)$ is a retraction of $\DGA$s then the associated diagram of vector spaces is injective. In particular, the minimal model of the cohomology diagram is level-wise minimal.
\end{prop}
\begin{proof}

Since the structure map $\cA_{e,G}:\cA(G/e)\to \cA(G/G)$ is a retraction, there exists $i:\cA(G/G)\to \cA(G/e)$ such that $\cA_{e,G}\circ i=id$. This  implies $\cA_{e,G}$ is surjective, and it follows that $\cA$ is injective diagram  of graded algebras. Note that for any minimal system of $\DGA$s $\sN$, the $\DGA$ $\sN(G/G)$ is non-equivariantly minimal by construction. 
    Let $\rho:\sM\to \cA$ be the minimal model and let $\sM_{e,G}:\sM(G/e)\to \sM(G/G)$ be the corresponding structure map. We claim that there exists $j:\sM(G/G)\to \sM(G/e)$ an inclusion map of $\DGA$s. 
    Since $\cA(G/e)$ is a $\DGA$ with zero differential, $\rho(G/e):\sM(G/e)\to \cA(G/e)$ is a surjective quasi-isomorphism, by the Lifting Lemma, there exists a lift $j:\sM(G/G)\to \sM(G/e)$ such that the diagram commutes  
        \[
\xymatrix{
\sM(G/e)\ar[r]^{ \rho(G/e)}  & \cA^i(G/e) \\
\sM(G/G) \ar[r]_{\rho(G/G)}\ar[u]^{j}   & \cA^i(G/G)\ar[u]_{i}
}
\]
Therefore, 
\begin{equation}
    \cA_{e,G}\circ i \circ \rho(G/G)= \cA_{e,G} \circ \rho(G/e)\circ j\\
    \implies \rho(G/G)=\rho(G/G)\circ \sM_{e,G}\circ j
\end{equation}
Since $\rho(G/G)$ is a quasi-isomorphism , $\sM_{e,G}\circ j:\sM(G/G)\to \sM(G/G)$ is a quasi-isomorphism. It then follows that, $\sM_{e,G}\circ j$ is an isomorphism and therefore $j:\sM(G/G)\to \sM(G/e)$ is an inclusion.

%Note that for any minimal system of $\DGA$s $\sN$, the $\DGA$ $\sN(G/G)$ is non-equivariantly minimal. Indeed, as $\sN$ is minimal it can be written as increasing union of elementary extensions $\sN_i =\mathbb{Q}(\oplus_{j=1}^{n_i} \underline{ V_j^i})$, where for each associated system of vector spaces $\underline{ V^i}$, we have the corresponding injective resolution $ \underline{ V_0^i}\to \underline{ V_1^i}\to\cdots\to  \underline{ V_{n_i}^i}$. For each $i$, $\underline{V_1^i}(G/G)=\underline{V}^i(G/G)$ from the definition of injective envelope. Therefore, the elementary extension's $G/G$-level is the same as that of the Hirsch extension (as in the non-equivariant case). 

Next, we show that all the associated systems of vector spaces are injective by induction on $n$ where $\sM=\cup \sM_m$. 
    Recall $\sM_n=\sM_{n-1}(\underline{V})$, where $\underline{V}$ is $H^{n+1}(\ker (\beta) \oplus \underline{\QQ})$ is the associated diagram of vector spaces at the $n$-th stage. Any element of $\underbar{V}(G/G)$ looks like the product of the elements of $\sM_{n-1}(G/G)$, $\cA(G/G)$ and $\sum \cA(G/G)$. We study case by case to conclude that $\underbar{V}(\hat{e}_{e,G})$ is surjective. Let $[x]\in \underbar{V}(G/G)$.

    \begin{enumerate}
        \item 
        If $x\in \sM^{n+1}_{n-1}(G/G)$ i.e., $[x]\in \underline{V}(G/G)$. Then $j(x)\in \sM_{n-1}(G/e)$. 

        As $\beta=\rho$ on $\sM_{n-1}$, we have $i\circ \rho(G/G)(x)=\rho(G/e)\circ j(x)$, this implies $i\circ \beta(G/G)(x)=\beta(G/e)\circ j(x)$, which implies $j(x)\in \ker \beta(G/e)$. 
        As $j$ is a $\DGA$-map we get $dj(x)=jd(x)$, which gives $j(x)\in \underline{V}(G/e)$. 
        \item If $x\in \sum \cA(G/G)$ then by injectivity of $\cA$, one gets a pre-image in $\sum \cA(G/e)$. The differential is zero for elements in $\sum \cA$. So we get a pre-image in $\underbar{V}(G/e)$. 

        \item Assume $x$ is the product of elements in $\sM_{n-1}$, $\cA$ and $\sum \cA$. In this case, note that the maps $i,j$ induce a $\DGA$-map $g: \sM(G/G)\otimes\QQ (\cA \oplus \sum \cA)(G/G)\to \sM(G/e)\otimes\QQ (\cA \oplus \sum \cA)(G/e)$. If $x=m.a.sb$ where $m\in \sM(G/G)$, $a\in \cA(G/G)$ and $sb\in \sum \cA(G/G)$, with $[x]\in \underline{V}(G/G)$, then we have $[g(m.a.sb)]\in \underline{V}(G/e)$.

    \end{enumerate}

Hence, $\underbar{V}$ is injective and it follows that  the minimal model $\sM$ is level-wise minimal. 
\end{proof}

\begin{exmp}\label{exmp:assosiated inj}

Consider $X=(S^3\vee S^3)\times S^5$  with action of $G$, where $G=C_2$ acts on $(S^3 \vee S^3)$ by switching copies and acts on $S^5$ trivially. We denote the cohomology diagram by $\cA$. 

Then $\cA(G/e)=\wedge(x,y)/<xy>\otimes \wedge (z)$ and $\cA(G/G)=\wedge (z)$ where  the generators of the cohomology are $x$, $y$ and $z$ and $deg(x)=3=deg(y)$ and $deg(z)=5$. 

Using Proposition \ref{prop:injcohom} we get that the given cohomology diagram is injective. Let $\sM$ denote its minimal model.

We have, 
\[
\sM_3=
\begin{cases}
     \wedge(a,b), & \text{at G/e}\\
    \mathbb{Q},  & \text{ at G/G}.
\end{cases}
\] 

Note that $\rho: \sM_3\to \cA$ is the extension of the map $\rho: \sM_2\to \cA$. It is defined by $\rho(a)=x$, $\rho(b)=y$.

It can be verified that  $\sM_4 \cong \sM_3$, as there are no elements in the cohomology of degree $4$.
 
 Computation for $\sM_5:$ From \cite{LS02}, the construction of minimal model we get
\[
\underbar{V}:=H^6(R)=
\begin{cases}
     \mathbb{Q}_6(sz,ab), & \text{at G/e}\\
    \mathbb{Q}_6(sz),  & \text{ at G/G},
\end{cases}
\] 

where $\QQ_n(x_i)$ denotes the $\QQ$-vector space generated by the elements $x_i$ of degree $n$. 

Using \Cref{prop:injcohom}, we see that $\underbar{V}$ is injective and can compute $\sM_5$.

One can show that, $$H^6(\sM_4;\underbar{V})\cong Hom(\underbar{V},H^6(\sM_4))=\mathbb{Q}_{\gamma},$$

where the map $\gamma: \underbar{V}\to H^6(\sM_4)$ takes $ab \to ab$ and $\QQ_{\gamma}\cong \mathbb{Q}$.

Thus, we have that,
\[
\sM_5=\sM_4^{\gamma}(\underbar{V})=
\begin{cases}
     \wedge(a_3,b_3,c_5,d_5), & \text{at G/e}\\
    \wedge(c_5),  & \text{ at G/G}.
\end{cases}
\] 

The differential is $d(a_3)=0=d(b_3),$ $d(d_5)=ab$ , $d(c_5)=0$, 
and $\rho :\sM_5\to \cA$ maps $c_5 \to z$ and $d_5\to 0.$

\end{exmp}

\section{Classifying rational homotopy types with isomorphic cohomology} \label{section:minimalmodels}

In this section, we consider the diagram of graded algebras,  $\cA^\ast$ over $\sO_{C_p}$, where $p$ is a prime number, and the structure map $\cA_{e,G}$ is a retract. Then by \Cref{Z_P injimpliesminimal}, we get that the minimal model for $\cA^{\ast}$ is level-wise minimal. We denote the minimal model for $\cA^{\ast}(G/H)$ by $N_H$. We use techniques from \cite{HIRO} to give an inductive construction for the minimal systems of $\DGA$s with cohomology algebra $\cA^{\ast}$. We assume  $\sM$ is the minimal model for $\cA^\ast$, so that $\sM(G/H)\cong N_H$, for each subgroup $H$ of $G$. We cannot directly extend the results of \cite{HIRO} because minimal algebras with cohomology $\cA^{\ast}(G/H)$, for every $H\leq G$, do not always give rise to a minimal system with cohomology diagram $\cA^{\ast}$.

The following proposition compares the level-wise weak equivalence of diagrams of $\DGA$s with a weak equivalence of the diagrams of $\DGA$s when they are level-wise minimal in the non-equivariant sense.
 \begin{prop}\label{prop:nonisomorphic}
Let $\sM$ and $\sN$ be two $\DGA$ diagrams over $\sO_G$, for some finite group $G$ having the same cohomology at each level of subgroups, and $\sM(G/H)$ and $\sN(G/H)$ are minimal algebras for every subgroup $H$ of $G$. Let $\sI_{\sM}$ and $\sI_{\sN}$ denote their injective envelopes. If $\sI_{\sM}$ and $\sI_{\sN}$ are weakly equivalent then for every subgroup $H$ of $G$, $\sM(G/H)$ and $\sN(G/H)$ have the same rational homotopy type.

Conversely, if there is a $\DGA$ diagram map $r:\sM\to \sN$ such that $r(G/H):\sM(G/H)\to \sN(G/H)$ is a weak equivalence for every subgroup $H$ of $G$, then $\sI_{\sM}$ and $\sI_{\sN}$ are weakly equivalent.
\end{prop}
\begin{proof}
Let $\sI_{\sM}$ and $\sI_{\sN}$ denote the injective envelopes for $\sM$ and, $\sN$ respectively. Let us assume that there is a weak equivalence $\phi:\sI_{\sM} \to \sI_{\sN}$. Then we have the diagram
 \[
\xymatrix{
\sM \ar[d]_{i}  & \sN \ar[d]^{j} \\
\sI_{\sM} \ar[r]  ^{\phi}    & \sI_{\sN}
}
\]
where $i,j$ are inclusions, which are also quasi-isomorphisms. Thus, we get a weak equivalence between $\sM$ and $\sN$ and for every subgroup $H$ of $G$, $\sM(G/H)$ and $\sN(G/H)$ are quasi-isomorphic as $\DGA$s. 

Conversely, if there exists a map $r:\sM\to \sN$ satisfying the hypothesis then consider the composition map $$\sM\xrightarrow{r} \sN\to \sI_{\sN}.$$ 
By \cite[Proposition 8]{FINE}, we get a map $\phi:\sI_{\sM}\to \sI_{\sN}$ which is the extension of the composition above.
\[
\xymatrix{
\sM\ar[r]^{r} \ar[d]_{i}  & \sN \ar[d]^{j} \\
\sI_{\sM} \ar@{-->}[r]_{\phi}   & \sI_{\sN}
}
\]

 Hence, this is a quasi-isomorphism, since other maps in the above diagram are quasi-isomorphism.
\end{proof}

We now describe the subset $\mathcal{N}^{C_p}_{\cA^{\ast}}$ of $\mathcal{M}^{C_p}_{\cA^\ast}$, which is defined to be the set of systems of $\DGA$s over $\sO_{C_p}$ which have cohomology $\cA^\ast$ up to weak equivalence and are minimal at each level $G/H$.\\

\begin{cons}\label{cons:plural}

Let $\cA^\ast$ be a diagram of graded algebras over $\sO_G$, where $G=C_p$ so that the structure map $\cA_{e,G}:\cA(G/e)\to \cA(G/G)$ is a retract.

We say that minimal algebras $p^i_{n-1}$, $i=1,2$ satisfy the Minimal ${(n-1)}$ property for $\cA^{\ast}$ if they satisfy the following conditions :  \\

\textbf{$(1)_{n-1}$: }Both the algebras $p^i_{n-1}$ for $i=1,2$, are generated by elements of degree $\leq n-1$ and for $i=1,2$ they are $\QQ[G/e]$ and $\QQ[G/G]$-modules respectively.

\textbf{$(2)_{n-1}$: }There is a $C_p$-$\DGA$-map $\delta: p^1_{n-1}\to p^2_{n-1}$ such that $\sP_{n-1}=\{p^1_{n-1},p^2_{n-1},\delta\}$ defines a diagram of $\DGA$s with  a morphism between diagram of graded algebras $$\sigma^j_{n-1}:(H^\ast(\sP_{n-1})(n))^j\to \cA^j$$

which is an isomorphism for $j\leq n-1$, and a monomorphism for $j=n$, where $((H^{\ast}\sP_{n-1})(n))$ is sub graded algebra diagram of $H^{\ast}(\sP_{n-1})$ generated by elements of degree $\leq n$.

 For instance, let $(\sM,\rho)$ be the minimal model of $\cA$, by \Cref{Z_P injimpliesminimal} it follows that $\sM(G/e)$ and $\sM(G/G)$ are minimal models for $\cA(G/e)$ and $\cA(G/G)$.  There are maps $\rho(G/e):\sM(G/e) \to \cA(G/e)$ and $\rho(G/G):\sM(G/G)\to \cA(G/G)$ which are quasi-isomorphisms.  Then the subalgebra $(\sM(G/H))_{n-1}$ of $\sM(G/H)$ generated by elements of degree less than or equal to $n-1$, satisfies the Minimal $(n-1)$ property for $\cA^{\ast}$  for $H$ equal to $\{e\}$ or $G$.

We denote by $\underline{K}$ the kernel of the map $\sigma^{\ast}_{n-1}:H^{\ast}((\sP_{n-1}(n))\to \cA^{\ast}$. Also, we let $K^1$ and $K^2$ denote $(\underline{K}(G/e))^{n+1}$ and $(\underline{K}(G/G))^{n+1}$ respectively.

%{\color{blue} Now, define $K_j$  to be the kernel of the maps $\sigma_{j,n-1}^{n+1}|H((p^j_{n-1})(n))^{n+1}$ and, for $j=1,2$. \textcolor{red}{The indices are wrong. maybe use $i$ for $\sigma$ ?} I want to remove this}

Let $p^i_D$ be the minimal algebras obtained by adding generators to $p^i_{n-1}$, whose differentials form a basis for $K^i$ for $i=1,2$. 

Consider the diagram of $\DGA$s $\sP_D=\{p^1_D,p^2_D,\theta\}$, where $\theta$ is induced from $\delta$ in $\sP_{n-1}$ and let 
$$\sigma_D:(H^{\ast}(\sP_D)(n))^{\ast}\to \cA^{\ast}$$ be the morphism induced by $\sigma_{n-1}$.

We set $$\dim_\QQ\cA(G/H)^{n+1}=u_i,~~\dim_{\QQ}\cA(G/H)^{n+1}/(\cA(G/H)(n))^{n+1}=s_i,$$ 
 for $(i,H)=(1,e)$ or $(2,G)$
 and for $i=1,2$ we set
 
 $$\dim_{\QQ}H^{n+1}(p^i_D)=v_i,~~\dim_{\QQ}\frac{H^{n+1}(p^i_D)}{(H^\ast(p^i_D)(n))^{n+1}}=t_i.$$
 
Since the sub-algebras generated by elements of degree $< n$ of the diagram of algebras $H^{n+1}(\sP_D)$ and $H^{n+1}(\sP_{n-1})$ are isomorphic, we have $$u_i-s_i=v_i-t_i~~for~~i=1,2.$$
 
 For $i=1,2$, let $l_i$ be integers satisfying $$\max(0,t_i-s_i)\leq l_i\leq t_i~~for~~i=1,2$$ and $W_i$ be subspaces of $H^{n+1}(p^i_D)$ such that $\dim W_i=l_i$ for $i=1,2$ with the condition that \beq \label{eq:222}
 W_i\cap (H^\ast(p^i_D)(n))^{n+1}=\{0\}.
 \eeq
Denote the minimal algebras obtained by adding $l_i$ generators (in degree $n$) to $p^i_D$, whose differentials span $W_i$'s by $p^{i,W_i}$.   

This implies that $H(p^{i,W_i})(n)=H(p^i_D)(n)$ for $i=1,2$.  We then have graded algebra maps $$\sigma_{i,D}:(H(p^{i,W_i})(n))\to \cA^\ast(G/H)$$ with
$$H^{n+1}(p^{i,W_i})\oplus W_i=H^{n+1}(p^i_D)~~~~for~i=1,2,$$ such that $$\dim_{\QQ}\frac{H^{n+1}(p^{i,W_i})}{(H(p^{i,W_i})(n))^{n+1}}=t_i-l_i\leq s_i=\dim_{\QQ}\frac{\cA^{n+1}(G/H)}{(\cA(n))^{n+1}(G/H)}$$

 for $(i,H)=(1,G)$ or $(2,e)$.

Let $p^{i,W_i}_n$ be the minimal algebras obtained by adding generators (in degree $n$) with zero differential to $p^{i,W_i}$ so that the cokernel of the map $$\sigma_{i,D}:(H(p^{i,W_i})(n))^n\to \cA^n(G/H)$$ becomes trivial for $(i,H)=(1,e)$ or $(2,G)$.
We get  graded algebra maps $$\sigma^{\ast}_{W_i,n-1}:(H(p^{i,W_i}_n)(n))^\ast\to \cA^\ast(G/H)$$
 such that for $i=1,2$, $\sigma^{\ast}_{W_j,n-1}$ are isomorphisms for $\ast\leq n$.

We say that the pair $(W_1,W_2)$ satisfies  Condition \textbf{(K)} if the following holds.

\begin{itemize}
    \item[(a)]

Given linear monomorphisms, 
$$
\phi_i:H^{n+1}(p^{i,W_i})/(H(p^{i,W_i})(n))^{n+1}\to \cA^{n+1}(G/H)/(\cA(G/H)(n))^{n+1} 
$$
where $(i,H)=(1,e)$ or $(2,G)$,
 the maps $\sigma_{W_i,n-1}\oplus \phi_i$ can be extended to graded algebra maps $$\sigma_{W_i,n}:(H(p^{i,W_i})(n+1))^\ast\to \cA^\ast(G/H),$$
 which are $\QQ(G/e)$ and $\QQ(G/G)$-module maps for $i=1,2$ respectively.

 \item[(b)] \label{conddd2}
There exists morphism of $\DGA$s, $\eta:p^{1,W_1}_n\to p^{2,W_2}_n$ such that $\sP_n^{(W_1,W_2)}=\{p_n^{1,W_1},p_n^{2,W_2},\eta\}$ is a diagram of $\DGA$s and the pair of maps $\{\sigma_{W_1,n},\sigma_{W_2,n}\}$ induce a morphism $$\sigma_{(W_1,W_2),n}:H^{\ast}(\sP_n^{(W_1,W_2)}(n+1))^{\ast}\to \cA^{\ast}$$

satisfying conditions $(1)_{n}$ and $(2)_{n}$.
\end{itemize}

\end{cons}

\begin{remark}

\begin{enumerate}

\item If there exists $W_i$'s such that $\dim_{\QQ}W_i=t_i$ the map $\eta$ always exists so that, condition $(1)_{n}$ and $(2)_{n}$ are satisfied. In fact, they will lead us to the minimal model for $\cA^\ast$. 

\item The pairs $(p^{1,W_1}_n,\sigma_{W_1,n})$ and $(p^{2,W_2}_n,\sigma_{W_2,n})$ depend on the sub-spaces $W_i$. If there are several choices of $W_i$, i.e., the choices of $l_i$'s are large, there is a higher possibility of finding a minimal system containing $p^1_{n-1}$ and $p^2_{n-1}$.

    \item Consider a diagram of graded algebras, $\cA$ which is equivariantly $(n-1)$ I.F. with minimal model $\sM$ such that $\sM_{n-1}(G/e)=p^1_{n-1}$ and $\sM_{n-1}(G/G)=p^2_{n-1}$.
   % In general, the existence of the map $\theta^\ast$ is a  strong condition. 
 The map $\delta$ always exists when, $l_i=t_i$ and gives the minimal model for $(H^{\ast}(\cA),0)$. If $\delta$ doesn't exist for any $l_i<t_i$ then what we can say is that there is no other minimal system with cohomology $\cA$, whose $n$-th stage at $G/H$, for both $H=e, G$  is a minimal algebra.

\end{enumerate}
\end{remark}

Let $Gr(v_i,l_i)(\QQ)$ be the set of rational points of the Grassmann manifolds of $l_i$-dimensional $\QQ$-subspaces in the $v_i$-dimensional spaces $H^{n+1}(p^1_D)$ and $H^{n+1}(p^2_D)$ respectively for $i=1,2$. Define
$$\sM_{(l_1,l_2)}:=\{(W_1,W_2)\in Gr(v_1,l_1)(\QQ)\times Gr(v_2,l_2)(\QQ)|~ W_i\cap (H^\ast(p^i_D)(n))^{n+1}=\{0\}\}$$
and $$\sO_{(l_1,l_2)}:=\{(W_1,W_2)\in \sM_{(l_1,l_2)}|~(W_1,W_2)~ \text{ satisfies Condition}~\textbf{(K)}\}.$$

 Let $G_1$ and $G_2$ be the group of $\DGA$ automorphisms of $p^1_D$ and $p^2_D$ respectively. Then $G_1$ and $G_2$ act on $H^{n+1}(p_D)$ and $H^{m+1}(p^1_D)$ respectively, and hence on $Gr(v_1,l_1)(\QQ)$ and $Gr(v_2,l_2)(\QQ)$. Let $(W_1,W_2)\in \sO_{l_1,l_2}$ and $(\phi_1,\phi_2)$ be an element of $G_1 \times G_2$.

Then it is easy to see that $(\phi_1,\phi_2)$ can be extended to $\DGA$ isomorphisms 
$$\phi_i:p^{i,W_i}_n\to p^{i,\phi_i(W_i)}_n$$ for $i=1,2$ respectively.

Conversely, let $(W_1\oplus W_2)$ and $(W_1'\oplus W_2')$ be $(l_1+l_2)$-dimensional subspaces of $H^{n+1}(p^1_D)\oplus H^{m+1}(p^2_D)$ such that there are $\DGA$ isomorphisms $$f_1:p^{1,W_1}_n \to p^{1,W_1'}_n$$ and $$f_2:p^{2,W_2}_n\to p^{2,W_2'}_n,$$ then $(f_1|p^1_D,f_2|p^2_D)\in G_1\times G_2$.  

However, it is not always true that such a pair will give rise to a morphism of $\DGA$ diagrams over $\sO_{C_p}$, that is, for such a pair there may not exist a $C_p$-$\DGA$ map $\theta$ such that Condition \textbf{K(b)} is satisfied. \\

Further, for isomorphic pairs $(W_1,W_2)$ and $(W'_1,W'_2)$ if there are two pairs of maps $\theta_1$ and $\theta_2$ such that diagram commutes
\[
\xymatrix{
p_n^{1,W_1}\ar[r]^{f_1} \ar[d]_{\theta_1}  & p_n^{1,W'_1} \ar[d]^{\theta_2} \\
p_n^{2,W_2} \ar[r]_{f_2}   & p_n^{2,W'_2}
}
\]
\begin{center}
    \label{figureeeeeee}
\end{center}
for a pair of automorphisms $(f_1,f_2)$, then the  $\DGA$ diagrams $X=(p_n^{1,W_1},p_n^{2,W_2},\theta_1)$ and $Y=(p_n^{1,W'_1},p_n^{2,W'_2},\theta_2)$ are weakly equivalent. In particular,  if we take their injective envelopes $\sI_X$ and $\sI_Y$, by Proposition \ref{prop:nonisomorphic} there is a map $\sI_X\to \sI_Y$ which is a weak equivalence.

Moreover, $p_n^{1,W_1}, p_n^{2,W_2}, p_n^{1,W'_1},p_n^{2,W'_2}$ are minimal algebras in the non-equivariant sense,  any quasi-isomorphism between them is an isomorphism.

 We now define an equivalence relation on the set $\sO_{l_1,l_2}$ as follows. Two elements $(W_1,W_2)$ and $(W'_1,W'_2)$ in $\sO_{l_1,l_2}$ are equivalent if there is a pair of $\DGA$ isomorphisms $(f_1,f_2)$ such that either one of the following diagrams commutes 
\[
\xymatrix{
p_n^{1,W_1}\ar[r]^{f_1} \ar[d]_{\theta_1}  & p_n^{1,W'_1} \ar[d]^{\theta_2} \\
p_n^{2,W_2} \ar[r]_{f_2}   & p_n^{2,W'_2}
}
\]
or 
\[
\xymatrix{
p_n^{1,W'_1}\ar[r]^{f_2} \ar[d]_{\theta_2}  & p_n^{1,W_1} \ar[d]^{\theta_1} \\
p_n^{2,W'_2} \ar[r]_{f_2}   & p_n^{2,W_2}
}
\]

We denote this relation by $\kappa$.\\

Using  Propositions, \ref{Z_P injimpliesminimal}, \ref{prop:nonisomorphic} and Construction \ref{cons:plural}, we have the following theorem
\begin{theorem}\label{theo:MAIN}

Let the diagram of graded algebras $\cA^\ast$ over $\sO_{C_p}$, where $p$ is a prime number, such that the structure map $\cA(G/e)\to \cA(G/G)$ is a retract. Let $\sM$ be the minimal model for $\cA^{\ast}$. The set of isomorphism classes of minimal algebras $\sN_k$ with the properties

\begin{enumerate}
    \item $\sN_k$ contains the minimal algebra $\sM_{k-1}$ satisfying $(1)_{k}$ and $(2)_{k}$, 
    \item $\sN_k(G/H)$ is minimal in the non-equivariant sense for every subgroup $H$ of $C_p$, 
\end{enumerate}

 corresponds bijectively to the disjoint union of orbit spaces $$X_{k}=\bigsqcup\sO_{l_1,l_2}/\kappa,$$

where $\kappa$ is the equivalent relation mentioned in Construction \ref{cons:plural}.

\end{theorem}

\begin{remark}

Let the diagram of graded algebras $\cA^\ast$ over $\sO_{C_p}$, where $p$ is a prime number, be such that the structure map is a retract. If $\cA^\ast$ is equivariantly $(n-1)$-intrinsically formal and $H^i(N_{C_p})=H^i(\cA)=0$ for $i \geq n$, then $H^i(p_n^{2,W_2})=0$ for $i\geq n$ and, one has Theorem \ref{thm:cardofM}.
\end{remark}

\begin{theorem}\label{thm:cardofM}
Let the diagram of graded algebras $\cA^\ast$ over $\sO_{C_p}$, where $p$ is a prime number, be such that the structure map $\cA(G/e)\to \cA(G/G)$ is a retract. Let $\sM$ be the minimal model for $\cA^{\ast}$. If $\cA^\ast$ is equivariantly $(n-1)$-intrinsically formal and $N_{C_p}^i=0$ for $i\geq n$, where $N_H$  is the minimal model for $\cA^{\ast}(C_p/H)$ for every subgroup $H$ of $C_p$, then the set of isomorphism classes of the minimal system $\sK_n$ containing the minimal system $\sM_{n-1}$ is determined by the set $\mathcal{M}_{\cA^\ast(C_p/e)}$.

Moreover, if $\cA^{\ast}$ is equivariantly $(n-1)$-I.F. and $\cA^{\ast}(C_p/e)=0$ for $j>n+1$. Then the cardinality of $\mathcal{M}_{\cA^{\ast}}^{C_p}$ is the same as the cardinality of $\mathcal{M}_{\cA^\ast(C_p/e)}$. 
\end{theorem}

\begin{proof}
By \Cref{Z_P injimpliesminimal} it follows that the minimal model of $\cA^{\ast}$ is level-wise minimal. Throughout the proof, we use $G $ to denote $C_p$.

 We first show that if $\sB$ is any system of $\DGA$s such that $H^{\ast}(\sB)=\cA^{\ast}$ and $\sB_{n-1}=\sM_{n-1}$ then $\sB$ is weakly equivalent to a system of $\DGA$s which is minimal at each $G/H$ in the non-equivariant sense.

To show this, let $N_H$ denote the non-equivariant minimal model for $\sM(G/H)$. Let $L$ denote the minimal model for $\sB(G/e)$. We define 
\[
\sV:=
\begin{cases}
     L, & \text{at G/e}\\
    N_G,  & \text{ at G/G}.
\end{cases}
\] 

We claim that $\sV$ is a system of $\DGA$s, and it is weakly equivalent to $\sB$.

First, note that $\sB(G/G)$ is quasi-isomorphic to $N_G=\sV(G/G)$. This is true since $\sB_{n-1}(G/G)$ contains $\sM_{n-1}(G/G)$ as a sub DGA in the non-equivariant sense. Also, $H^k(\sB(G/G))=H^k(N_G)=\cA(G/G)$ and $H^i(\sB(G/G))=0$ for $i\geq n$. Thus, $\sB(G/G)$ is obtained by adding an acyclic vector space to $N_G$, and hence the inclusion map $N_G \to \sB(G/G)$, is a quasi-isomorphism. Now observe that, $(\sV(G/e))^{n-1}=L^{n-1}=N_e^{n-1}$ and $N_G^i=0$ for $i\geq n,$ so we can define the map $\theta: \sV(G/e) \to \sV(G/G)$ by $\theta|\sV(G/e)^{n-1}$ as the map $\sM_{n-1}(G/e)=N_e^{n-1}\to \sM_{n-1}(G/G)= N_G^{n-1}$ and zero at $(\sV(G/e)^{i})$ for $i\geq n$, so that $\sV$ is a diagram $\DGA$s.

Next, we show that $\sV$ is a system of $\DGA$s. Note that $\sV_{n-1}=\sM_{n-1}$ is already a system of $\DGA$s. Since $\sV^i(G/G)=N^i_G=0$ for $i\geq n$, so the maps $\sV^i(G/e) \to \sV^i(G/G)$ are zero maps for $i\geq n$ and hence surjective. $\sV_{n-1}=\sM_{n-1}$ being a system, the map $\sV^j(G/e)\to \sV^j(G/G)$ is surjective for $j\leq n-1$ thus, combining we get that $L \to N_G$ is surjective and hence, by Proposition \ref{prop:injcohom}, $\sV$ is a system.

Also, we have a map from $\sV \to \sB$ which is identity up-to the $(n-1)$th stage and since $\sV^i(G/G)=0$ for $i\geq n$ so $\sV^i(G/e) \to \sB^i(G/e)$ the quasi-isomorphisms will fit in to the diagram $\sV \to \sB$. 
\[
\xymatrix{
\sV^i(G/e)\ar[r] \ar[d]_{\sV(\hat{e}_{e,G})}  & \sB^i(G/e) \ar[d]^{\sB(\hat{e}_{e,G})} \\
\sV^i(G/G) \ar[r]   & \sB^i(G/G)
}
\]

Thus, using the converse of Proposition \ref{prop:nonisomorphic} we get that $\sV$ and $\sB$ are weakly equivalent.

Thus, any such $\DGA$ diagram $\sK$ of the form 
\[
\sK:=
\begin{cases}
     m, & \text{at G/e}\\
    N_G,  & \text{ at G/G}.
\end{cases}
\] 

with the property that $m$ is a minimal algebra, and $\sK_{n-1}=\sM_{n-1}$ is a system of $\DGA$s. Thus, the plural homotopy types are determined by $m$. The cardinality of the set of such $m$ is determined by the set $\cM_{\cA^{\ast}(G/e)}$.

\end{proof}

\begin{corollary}\label{cor:main}

Let the diagram of graded algebras $\cA^\ast$ over $\sO_{C_p}$, where $p$ is a prime number, be such that its structure map is a retract. Let $\sM$ be the minimal model for $\cA^{\ast}$. If $\cA^\ast$ is equivariantly $(n-1)$-intrinsically formal, $N_{C_p}^i=0$ for $i\geq n$ and $N_e$ is formal, then $\cA^{\ast}$ is equivariantly intrinsically formal. In particular, $\sM$ is equivariantly formal.  
\end{corollary}
\begin{proof}
By Theorem \ref{thm:cardofM} $\sM_{n-1}$ is determined by the set $\cM_{\cA^{\ast}(C_p/e)}$ and since $H^{\ast}(N_e)$ is formal, the cardinality of this set in one. This holds for each $n$ and hence $\cA^{\ast}$ is equivariantly I.F. 
\end{proof}

\section{Examples}\label{sec:examples}

\begin{exmp} \label{example sec5}
 Let $n\geq 2$  and $p$ be a prime. Let $X= S^n \times \cdots \times S^n$ be the product of $p$ copies of $S^n$ with the  $G=C_p$ action $t(a_1,a_2,\cdots, a_p)= (a_2,a_3,\cdots, a_{p-1},a_1)$ where $t$ is a generator of $C_p$. Then we have the following two cases.
\begin{enumerate}

\item\label{oddddd}

 Let $n $ be odd. From Example \ref{odd:product} we conclude that the cohomology diagram is equivariantly $n$-I.F. Using Corollary \ref{cor:main} we conclude that the product of odd spheres under the above action is  equivariantly formal.

\item\label{evennnn} Let $n $ be an even integer. From Example \ref{even:product} the cohomology diagram is equivariantly $(2n-1)$ I.F. So by using Corollary \ref{cor:main} we conclude that the product of even spheres under the above action is equivariantly formal.
\end{enumerate}  
\end{exmp}
    \begin{exmp}
        
    \label{oddwedgecross} Consider $X=(S^3\vee S^3)\times S^5$ with diagonal action of $G=C_2$, where $C_2$ acts on $(S^3 \vee S^3)$ by switching copies and acts trivially on $S^5$.
     
We denote the cohomology diagram by $\cA$, the generators of the cohomology by $x$, $y$, and $z$. Note that, $deg(x)=3=deg(y)$ and $deg(z)=5$, and
\[
\cA:=
\begin{cases}
    \wedge(x,y)/<xy>\otimes \wedge (z), & \text{at G/e}\\
    \wedge (z),  & \text{ at G/G}.
\end{cases}
\] 

First, note that $\cA$ is injective by \Cref{prop:injcohom}. We claim that the cohomology diagram is equivariantly $6$-I.F. 

To see this, first, note that the cohomology diagram is equivariantly $3$-I.F. since $\cA$ is $2$-connected and the minimal system of $\DGA$s 
 \[
\sM_3=
\begin{cases}
     \wedge(a_3,b_3), & \text{at G/e}\\
    \QQ,  & \text{ at G/G}
\end{cases}
\] 
which is equivariantly I.F. 

Next, if we consider any minimal system $\sM'$ with $H^{\ast}(\sM')=\cA$ then $\sM'_3=\sM_3$. Now $\sM'_4=\sM_3(\underbar{V})$.  

We claim that the inclusion $\sM_3 \to \sM'_4$ is a $4$-isomorphism, so $\cA$ becomes equivariantly $4$-I.F. 

Since $\sM'_4\to \sM'$ is $4$-isomorphism, the map $H^i(\sM'_4)\to H^i(\sM')\cong H^i(\sM)=\cA^i$ is isomorphism for $i\leq 4$ and monomorphism for $i=5$. 

Therefore, $H^i(\sM_3) \to H^i(\sM'_4)$ is isomorphism for $i\leq 4$ as $H^4(\sM')=0$. Hence, we need to show $H^5(\sM_3)\to H^5(\sM'_4)$ is monomorphism, i.e., we must show $H^5(\sM'_4)=0$. 

Assume to the contrary that there exists $[\beta]\in H^5(\sM'_4)$ non-zero. Then $\beta \in (\sM_3(\underbar{V}))^5$. Now we consider the corresponding resolution of $\underbar{V}$
\[
\xymatrix{
\underline{V}\ar[d]_{\alpha} \ar[r] & \underline{V}_0\ar[d]^{ \alpha_0}\ar[r]^{w_0}& \underline{V}_1\ar[d]^{\alpha_1}\ar[r]^{w_1} & \underline{V}_2\ar[r]\ar[d]^{\alpha_2} & \cdots  \\
\underline{Z}^{4}(\sM_3)\ar[r]         & \sM_3^4\ar[r]_{d}&\sM_3^5\ar[r]_d & \sM_3^6 \ar[r] &\cdots
}
\] 

Note that $\underbar{Z}^4(\sM_3)=\sM^4_3=\sM^5_3=0$. Thus, $\alpha_0=\alpha_1=0$ and the differential is dependent only on $w_0,w_1$ at the $\underbar{V}_0$ and $\underbar{V}_1$ level, respectively.
The element $\beta \in \underbar{V}_1$ is non-zero in the cohomology so it is not in the image of $w_0$. In particular, $\beta \in \underline{V}_1/\Im (w_0)$ and $ \underline{V}_1/\Im (w_0)\to\underline{V}_2$ is an injection, implies that $d(\beta)\neq 0$ and $\beta$ is not a cohomology class. This implies $H^5(\sM'_4)=0$. 
Hence, $H^5(\sM_3)\to H^5(\sM'_4)$ is a monomorphism. So the map $\sM_3 \to \sM'_4$ is $4$-isomorphism and hence $\cA$ is equivariantly $4$-I.F. 

Next, we claim that $\cA$ is equivariantly $5$-I.F. 

We note that 
\[
\sM_5=
\begin{cases}
     \wedge(a_3,b_3,c_5,d_5), & \text{at G/e}\\
    \wedge(c_5),  & \text{ at G/G}
\end{cases}
\] 

So if $\sM'$ is any minimal system with $H^{\ast}(\sM')=\cA$, then by previous argument we get $\sM'_3=\sM_3$ and $\sM_4 \to \sM'_4$ is a $4$-isomorphism.

We claim that $\sM_5 \to \sM'_5$ is a $5$-isomorphism.  

Since $\sM'_5 \to \sM'$ is a $5$-isomorphism, there are elements $p,q \in \sM'_5$ of degree $5$ so that $[p]$ maps to $z$ and the differential of $q$ kills the product $xy$ and $H^6(\sM'_5)=0$.

We define the map $\sM_5 \to \sM'_5$ which when restricted to $\sM_3$ is the previous map and sends $a,b$ to $p,q$ respectively. The induced map indeed is an isomorphism $H^i(\sM_5)\to H^i(\sM'_5)$ for $i\leq 5$ and monomorphism for $i=6$.

Hence, $\cA$ is equivariantly $5$-I.F. Since $\cA^6= \cA^7 =0$ from Proposition \ref{prop:IIF} we conclude that $\cA$ is equivariantly $6$-I.F.

The minimal system up to $6$-th stage is given by 
 \[
\sM_6=
\begin{cases}
     \wedge(a_3,b_3,c_5,d_5), & \text{at G/e}\\
    \wedge(c_5),  & \text{ at G/G}
\end{cases}
\] 

The differential is $d(a_3)=0=d(b_3),$ $d(d_5)=ab$ , $d(c_5)=0$, 
and $\rho :\sM_6\to \cA$ maps $c_5 \to z$ and $d_5\to 0.$

 Thus, the cohomology diagram satisfies the hypothesis of Theorem \ref{thm:cardofM}. Hence, the set of isomorphism classes of the minimal system containing $\sM_6$ is in bijection with $\cM_{\cA^{\ast}(G/e)}$. From \cite{HIRO} and Theorem \ref{thm:cardofM}, it follows that the $\cM_{\cA^{\ast}}^{C_2}$ consists of three points. The description is given in \cite[Example 5]{HIRO}. 
 \end{exmp}
 \begin{exmp}\label{exmp: 3}
 
 Consider the space $X=(S^2\vee S^2)\times S^3$ with the action of $G=C_2$, by switching copies of $(S^3 \vee S^3)$ and acts on $S^5$ trivially.
We denote the cohomology diagram by $\cA$, and the generators of the cohomology by $x$, $y$, and $z$. Note that, $deg(x)=2=deg(y)$ and $deg(z)=3$.
The cohomology diagram is given by
\[
\cA=
\begin{cases}
    \wedge(x,y)/<xy>\otimes \wedge (z), & \text{at G/e}\\
    \wedge (z),  & \text{ at G/G}.
\end{cases}
\] 

By \Cref{prop:injcohom}, we conclude that $\cA$ is injective.
The cohomology diagram is equivariantly $3$-I.F. The minimal system up to $3$-rd stage is given by 
 \[
\sM_3=
\begin{cases}
     \wedge(a_2,b_2,c_3,d_3,e_a,e_b), & \text{at G/e}\\
    \wedge(c_3),  & \text{ at G/G}
\end{cases}
\] 

The differential is given by $d(a_2)=0=d(b_2),$ $d(d_3)=ab$ , $d(c_3)=0$, $d(e_a)=a^2$, $d(e_b)=b^2$
and $\rho :\sM_3\to \cA$ maps $(c_3,a_2,b_2) \to (z,x,y)$ and $(d_3,e_a,e_b)\to (0,0,0).$

 Thus, the cohomology diagram satisfies the hypothesis of \Cref{thm:cardofM}. Hence, the set of isomorphism classes of the minimal system containing $\sM_4$ is determined by the set $\cM_{\cA^{\ast}(G/e)}$. From \Cref{thm:cardofM}, we conclude that $\cM_{\cA^{\ast}}^{C_2}$ is the same as $\cM_{\cA^{\ast}(G/e)}$. By \cite[Example 4]{HIRO} then we conclude that $\cM_{\cA^{\ast}}^{C_2}$ contains exactly two points.
 
 \end{exmp}

\begin{exmp}\label{wedgeeeodd} There is a diagonal $G$ action on $X=S^3\vee S^3 \vee S^5$ where $G=C_2$ acts on $(S^3 \vee S^3)$ by switching copies and acts on $S^5$ trivially.

Note that the fixed point set $X^G$ is homeomorphic to $S^5$.
We denote the cohomology diagram by $\cA$, and the generators of the cohomology by $x$, $y$, and $z$. Note that, $deg(x)=3=deg(y)$ and $deg(z)=5$, and

\[
\cA=
\begin{cases}
     \mathbb{Q}(x,y,z)/<x^2, z^2, xy,xz,yz>, & \text{at G/e}\\
    \mathbb{Q}(z)/<z^2>,  & \text{ at G/G}
\end{cases}
\] 

This diagram of graded algebras is equivariantly $5$-I.F.
\[
\sM_5=
\begin{cases}
     \wedge(a_3,b_3,c_5,d_5), & \text{at G/e}\\
    \wedge(c_5),  & \text{ at G/G}
\end{cases}
\] 
The differential is $d(a_3)=0=d(b_3),$ $d(d_5)=ab$ , $d(c_5)=0$,

Thus, the cohomology diagram satisfies all the hypotheses of Theorem \ref{thm:cardofM} and hence up to isomorphism the number of minimal algebras containing $\sM_5$ with the same cohomology diagram can be computed similarly. 
\end{exmp}
\begin{exmp}\label{exmp:last}

Let $n>1$ be an integer and let $X=(S^n\vee S^n)\times S^{2n-1}$. Then there is a $G=C_2$ action on $X$ switching the two copies of $S^n$ and keeping the $S^{2n-1}$ copy fixed. 
We denote the cohomology diagram by $\cA$ and the generators of the cohomology by $x$, $y$, and $z$. Note that, $deg(x)=n=deg(y)$ and $deg(z)=2n-1$, and
\[
\cA:=
\begin{cases}
    \wedge(x,y)/<xy>\otimes \wedge (z), & \text{at G/e}\\
    \wedge (z),  & \text{ at G/G}.
\end{cases}
\] 

Since the map $\cA(G/e) \to \cA(G/G)$ is onto, we see that the given cohomology diagram is injective. 

\begin{enumerate}
    \item If $n$ is odd, the cohomology diagram is equivariantly $2n$-intrinsically formal. The minimal system up to $2n$ stage is given by

     \[
\sM_{2n}=
\begin{cases}
     \wedge(a_n,b_n,c_{2n-1},d_{2n-1}), & \text{at G/e}\\
    \wedge(c_{2n-1}),  & \text{ at G/G}
\end{cases}
\] 

The differential is $d(a_n)=0=d(b_n),$ $d(d_{2n-1})=ab$ , $d(c_{2n-1})=0$, 
and $\rho :\sM_{2n}\to \cA$ maps $c_{2n-1} \to z$ and $d_{2n-1}\to 0.$

 Thus, the cohomology diagram satisfies the hypothesis of Theorem \ref{thm:cardofM}. Hence, the set of isomorphism classes of the minimal system containing $\sM_6$ is obtained similarly.

 \item If $n$ is even, the cohomology diagram is equivariantly $(2n-1)$-I.F. The minimal system up to the $(2n-1)$-th stage is given by 

 \[
\sM_{2n-1}=
\begin{cases}
     \wedge(a_n,b_n,c_{2n-1},d_{2n-1},e_a,e_b), & \text{at G/e}\\
    \wedge(c_{2n-1}),  & \text{ at G/G}
\end{cases}
\] 

The differential is $d(a_n)=0=d(b_n),$ $d(d_{2n-1})=ab$ , $d(c_{2n-1})=0$, $d(e_a)=a^2$, $d(e_b)=b^2$
and $\rho :\sM_{2n-1}\to \cA$ maps $(c_{2n-1},a_n,b_n) \to (z,x,y)$ and $(d_{2n-1},e_a,e_b)\to (0,0,0).$

 Thus, the cohomology diagram satisfies the hypothesis of Theorem \ref{thm:cardofM}. Hence, the set of isomorphism classes of the minimal system containing $\sM_{2n}$ is in bijection with the set $\cM_{\cA^{\ast}(G/e)}$. The case $n=2$ is computed in \cite{HIRO}. 
\end{enumerate}

\end{exmp}

 \section*{Acknowledgements:} The authors would like to thank the referee and editor for their comments and suggestions, which helped improve the paper exposition. The first author would like to thank Samik Basu for sharing his insights on the subject.
 The second author would like to thank Aleksandar Milivojevic for some very helpful discussions on the proof of  Proposition \ref{prop:kiso implies quasi iso}. The second author was partially funded by U.G.C., India, through Senior Research Fellowship and IRCC fellowship (IIT Bombay) during this work.

\end{document}